\documentclass[12pt]{article}
\usepackage[mathscr]{eucal}
\begin{document}
\newtheorem{prop}{}[section]
\newtheorem{defi}[prop]{}
\newtheorem{lemma}[prop]{}
\newtheorem{rema}[prop]{}
\newcommand{\boma}[1]{{\mbox{\boldmath $#1$} }}
\def\J{{\mathscr J}}
\def\H{{\mathscr H}}
\def\K{{\mathscr K}}
\def\Kp{{\mathscr K}'}
\def\kp{k'}
\def\scrscr{\scriptscriptstyle}
\def\scr{\scriptstyle}
\def\dd{\displaystyle}
\def\z{\overline{z}}
\def\B{ B_{\mbox{\scriptsize{\textbf{C}}}} }
\def\Bc{ \overline{B}_{\mbox{\scriptsize{\textbf{C}}}} }
\def\ppartial{\overline{\partial}}
\def\d{\hat{d}}
\def\TT{T}
\def\G{ {\textbf G} }
\def\Hinf{ H^{\infty}(\reali^d, \complessi) }
\def\Hn{ H^{n}(\reali^d, \complessi) }
\def\Hm{ H^{m}(\reali^d, \complessi) }
\def\Ha{ H^{\d}(\reali^d, \complessi) }
\def\Ld{L^{2}(\reali^d, \complessi)}
\def\Lpi{L^{p}(\reali^d, \complessi)}
\def\Lq{L^{q}(\reali^d, \complessi)}
\def\Lr{L^{r}(\reali^d, \complessi)}
\def\Knb{K^{best}_n}
\def\k{\mbox{{\tt k}}}
\def\x{\mbox{{\tt x}}}
\def\D{\mbox{{\tt D}}}
\def\g{ {\textbf g} }
\def\QQQ{ {\textbf Q} }
\def\AAA{ {\textbf A} }
\def\gr{\mbox{graph}~}
\def\Q{$\mbox{Q}_a$~}
\def\PZ{$\mbox{P}^{0}_a$~}
\def\PZAL{$\mbox{P}^{0}_\alpha$~}
\def\PL{$\mbox{P}^{1/2}_a$~}
\def\PU{$\mbox{P}^{1}_a$~}
\def\PK{$\mbox{P}^{k}_a$~}
\def\PKU{$\mbox{P}^{k+1}_a$~}
\def\PI{$\mbox{P}^{i}_a$~}
\def\Pell{$\mbox{P}^{\ell}_a$~}
\def\PTM{$\mbox{P}^{3/2}_a$~}
\def\AZ{$\mbox{A}^{0}_r$~}
\def\AU{$\mbox{A}^{1}$~}
\def\epsilona{\epsilon^{\scriptscriptstyle{<}}}
\def\epsilonb{\epsilon^{\scriptscriptstyle{>}}}
\def\lgraffa{ \mbox{\Large $\{$ } \hskip -0.2cm}
\def\rgraffa{ \mbox{\Large $\}$ } }
\def\restriction{ \stackrel{\setminus}{~}\!\!\!\!|~}
\def\M{{\scriptscriptstyle{M}}}
\def\m{m}
\def\Fre{Fr\'echet~}
\def\I{{\mathcal N}}
\def\ap{{\scriptscriptstyle{ap}}}
\def\fiap{\varphi_{\ap}}
\def\BBB{ {\textbf B} }
\def\EEE{ {\textbf E} }
\def\FFF{ {\textbf F} }
\def\TTT{ {\textbf T} }
\def\KKK{ {\textbf K} }
\def\FFi{ {\bf \Phi} }
\def\GGam{ {\bf \Gamma} }
\def\a{a}
\def\ep{\epsilon_d}
\def\parn{\par\noindent}
\def\teta{M}
\def\elle{L}
\def\ro{\rho}
\def\al{\alpha}
\def\si{\sigma}
\def\be{\beta}
\def\ga{\gamma}
\def\de{\delta}
\def\la{\landa}
\def\te{\vartheta}
\def\ch{\chi}
\def\et{\eta}
\def\complessi{{\textbf C}}
\def\reali{{\textbf R}}
\def\interi{{\textbf Z}}
\def\naturali{{\textbf N}}
\def\T{{\textbf T}}
\def\T1{{\textbf T}^{1}}
\def\EE{{\mathcal E}}
\def\FF{{\mathcal F}}
\def\GG{{\mathcal G}}
\def\PP{{\mathcal P}}
\def\QQ{{\mathcal Q}}
\def\Np{{\hat{N}}}
\def\Lp{{\hat{L}}}
\def\Jp{{\hat{J}}}
\def\Pp{{\hat{P}}}
\def\Pip{{\hat{\Pi}}}
\def\Vp{{\hat{V}}}
\def\Ep{{\hat{E}}}
\def\Fp{{\hat{F}}}
\def\Gp{{\hat{G}}}
\def\Ip{{\hat{I}}}
\def\Tp{{\hat{T}}}
\def\Mp{{\hat{M}}}
\def\La{\Lambda}
\def\Ga{\Gamma}
\def\Si{\Sigma}
\def\Upsi{\Upsilon}
\def\Gag{{\check{\Gamma}}}
\def\Lap{{\hat{\Lambda}}}
\def\Sip{{\hat{\Sigma}}}
\def\Upsig{{\check{\Upsilon}}}
\def\Kg{{\check{K}}}
\def\ellp{{\hat{\ell}}}
\def\j{j}
\def\jp{{\hat{j}}}
\def\Stir{{\mathscr S}}
\def\BB{{\mathcal B}}
\def\LL{{\mathcal L}}
\def\SS{{\mathcal S}}
\def\DD{{\mathcal D}}
\def\VV{{\mathcal V}}
\def\WW{{\mathcal W}}
\def\OO{{\mathcal O}}
\def\CC{{\mathcal C}}
\def\RR{{\mathcal R}}
\def\AA{{\mathcal A}}
\def\CC{{\mathcal C}}
\def\JJ{{\mathcal J}}
\def\NN{{\mathcal N}}
\def\WW{{\mathcal W}}
\def\HH{{\mathcal H}}
\def\XX{{\mathcal X}}
\def\YY{{\mathcal Y}}
\def\ZZ{{\mathcal Z}}
\def\UU{{\mathcal U}}
\def\CC{{\mathcal C}}
\def\XX{{\mathcal X}}
\def\RR{{\mathcal R}}
\def\cir{{\scriptscriptstyle \circ}}
\def\circa{\thickapprox}
\def\vain{\rightarrow}
\def\ss{s}
\def\vains{\stackrel{\ss}{\rightarrow}}
\def\parn{\par \noindent}
\def\salto{\vskip 0.2truecm \noindent}
\def\spazio{\vskip 0.5truecm \noindent}
\def\vs1{\vskip 1cm \noindent}
\def\fine{\hfill $\diamond$ \vskip 0.2cm \noindent}
\newcommand{\rref}[1]{(\ref{#1})}
\def\beq{\begin{equation}}
\def\feq{\end{equation}}
\def\beqq{\begin{eqnarray}}
\def\feqq{\end{eqnarray}}
\def\barray{\begin{array}}
\def\farray{\end{array}}
\makeatletter
\@addtoreset{equation}{section}
\renewcommand{\theequation}{\thesection.\arabic{equation}}
\makeatother
\begin{titlepage}
\begin{center}
{\huge Quantitative functional calculus in Sobolev
spaces.}
\end{center}
\vspace{1truecm}
\begin{center}
{\large
Carlo Morosi${}^1$, Livio Pizzocchero${}^2$} \\
\vspace{0.5truecm}
${}^1$ Dipartimento di Matematica, Politecnico di
Milano, \\ P.za L. da Vinci 32, I-20133 Milano, Italy \\
e--mail: carmor@mate.polimi.it \\
${}^2$ Dipartimento di Matematica, Universit\`a di Milano\\
Via C. Saldini 50, I-20133 Milano, Italy\\
and Istituto Nazionale di Fisica Nucleare, Sezione di Milano, Italy \\
e--mail: livio.pizzocchero@mat.unimi.it
\end{center}
\vspace{1truecm}
\begin{abstract} In the
framework of Sobolev (Bessel potential) spaces
$H^n(\reali^d, \reali~ \mbox{or}~\complessi)$,
we consider the nonlinear Nemytskij operator
sending a function
$x \in \reali^d \mapsto f(x)$
into a
composite function $x \in \reali^d \mapsto G(f(x), x)$.
Assuming sufficient smoothness for $G$, we give a
"tame" bound on the $H^n$ norm of this composite function in terms
of a linear function of the $H^n$ norm of $f$, with a coefficient
depending on $G$ and on the $H^a$ norm of $f$, for all integers
$n, a, d$ with $a > d/2$. In comparison
with previous results on this subject, our bound is fully explicit,
allowing to estimate quantitatively the $H^n$ norm of
the function $x \mapsto G(f(x),x)$. When applied to the case $G(f(x), x) =
f^2(x)$, this bound agrees with a previous result of ours on the pointwise
product of functions in Sobolev spaces.
\end{abstract}
\vspace{1truecm} \noindent \textbf{Keywords:} Sobolev spaces,
inequalities, Nemytskij operators. \par \vspace{0.4truecm}
 \par
\vspace{2truecm} \noindent October 2001 (with updated references).
\par \noindent
\end{titlepage}
\section{Introduction.}
\label{intro}
Functional calculus in Sobolev spaces has been extensively discussed, from
the classic papers \cite{Ada2} \cite{Bou} to the recent book \cite{Run}.
These and the other works of our knowledge have been
mainly devoted to the well definedness of nonlinear
composition operators (Nemytskij operators),
and to the qualitative structure of the norm estimates for these operators.
On the contrary, in this paper we are interested in giving
fully quantitative norm estimates for the Nemytskij operators,
for the special case of $L^2$ type Sobolev spaces, including
evaluation of all numerical constants involved therein;
the spirit is the same as in the previous work \cite{mp2} on
the pointwise product of two
functions in a Sobolev space, whose results imply
in an elementary way estimates on the operators of
composition with polynomials. Here we will infer estimates for
much more general Nemytskij operators, covering essentially
the composition with any sufficiently derivable mapping. When the present,
general results are applied
to polynomial mappings, in the quadratic
case they reproduce the upper bounds
arising from \cite{mp2} (whose reliability was
also tested in the cited work, by comparison with the lower bounds obtained
from appropriate trial functions);
for polynomial Nemytskij operators of degree $\geq 3$, the upper bounds
derived in the present framework
are even more efficient than the elementary iteration
of the estimates in \cite{mp2}.
\parn
The Sobolev spaces we consider in this paper are the $L^2$ type Sobolev
(Bessel potential) spaces $H^{n}(\reali^d, \complessi)$,
or their real analogues $H^n(\reali^d, \reali)$,
with the corresponding norms $\|~\|_{n} := \|~\|_{H^n}$
(see \cite{Smi} \cite{Maz} or Sect.\ref{desc}); the attention is restricted to
the case of integer order $n$.
\parn
Before describing our results on
Nemytskij operators, let us fix some notations.
Consider a function $G : \BB \times \reali^d \vain \reali~ \mbox{or}~
\complessi$,
where $\BB$ is an interval of $\reali$
or a ball of $\complessi$, say, with center $0$; if
$f : \reali^d \mapsto \BB$ and $\chi : \reali^d
\mapsto \reali^d$, one can define the composition
$G(f,\chi) : \reali^d \mapsto \reali~\mbox{or}~\complessi$,
$x \mapsto G(f(x),\chi(x))$. We choose for $\chi$ the identity
mapping $\x : \reali^d \vain \reali^d$, $x = (x_1, ..., x_d)
\mapsto x$; thus
\beq G(f,\x) : \reali^d \mapsto \reali~\mbox{or}~
\complessi~, \qquad x \mapsto G(f(x), x)~. \label{defnem} \feq
The (generally nonlinear) mapping
\beq f \mapsto G(f, \x) \feq
will be called the Nemytskij operator associated to $G$. Its domain
is the set of functions (modulo equality almost everywhere)
$f : \reali^d \vain \BB$, $x \mapsto f(x)$, and the codomain
is made of functions $\reali^d \mapsto \reali~\mbox{or}~ \complessi$
(in the sequel, the expression "$G(f,\x)$ is defined " will always be
employed to mean that a function $f$ on $\reali^d$ takes values in
$\BB$). \parn
The classical problem of functional calculus
in Sobolev spaces is to prove that $G(f,\x)$ is in
a Sobolev space, and to estimate its norm,
when $G$ is sufficiently smooth and $f$ is in a related Sobolev space.
To discuss this topic one can fix a base integer $a > d/2$, so that
functions in $H^a$ be bounded and continuous (in applications,
$a$ is often chosen to be the smallest integer $> d/2$); a second
integer $n$ runs freely over $\naturali$,
labelling the scale of all $H^n$ spaces. It is known that
if $G$ satisfies suitable conditions of smoothness and boundedness,
$G(0,\x)$ is in $H^n$ and $f$ is in $H^n \cap$ a ball of $H^a$, then $G(f,\x)$
is defined, it belongs to $H^n$ and
\beq \| G(f,\x) - G(0,\x) \|_n \leq \Upsilon_{n a d}(G, \| f \|_a) \| f \|_n~;
\label{tame} \feq
as it is seen, the r.h.s. in this inequality is a linear
function of $\| f \|_n$, with a coefficient depending on $\| f \|_a$
in a way determined by $n, a, d$ and by the function $G$. \parn
Our aim is to provide an \textsl{explicit} expression for the function
$\Upsilon_{n a d}$ ruling this dependence:
this will be given in Prop.s \ref{nemreale} and \ref{nemc},
which are the main results of the paper.
One feature
of the expression we will obtain is that it does \textsl{not} depend dramatically
on the dimension $d$ (differently from the bounds frequently appearing
in similar situations, when the indices of derivatives are manipulated too
naively); this result is obtained by systematic use of tensorial methods
to express the derivatives of any order of functions on $\reali^d$ and
their composition with $G$.
Another feature of our result, already stressed, is that
it contains, as a subcase, the outcomes of \cite{mp2}
on the case $G(f, \x) = f^2 =$ the pointwise product of
$f$ by itself (estimates on pointwise product were derived in
\cite{mp2} using the Fourier transform, instead of the present tensor
approach, as an alternative tool to deal efficiently with the $d$-dependence).
\parn
In Eq.\rref{tame} one recognises the
structure of the so-called "tame" estimates, forming the basis of the
Nash-Moser implicit function theorem \cite{Ham}.
More precisely, we have a "tame estimate
of order zero" (an estimate of general order $r$ would
give the $n$-th norm of a nonlinear operator on $f$ in terms
of a linear function of $\| f \|_{n+r}$, with a coefficient
depending on $\| f\|_{a}$ for some fixed order $a$). Such zero order
estimates are used in \cite{quad} in relation to semilinear
evolution equations on the intersection
$\cap_{n=0}^{\infty} H^n$; the results of the
present paper allow quantitative estimates on the existence times
for the solutions, when the nonlinear part of the evolution equation
is a Nemytskij operator. \parn
Let us illustrate the organisation of the paper. In Sect.\ref{desc}
we describe the main results, i.e., the expression of $\Upsilon_{n a d}$
in the case of real and complex Sobolev spaces; this requires
some preliminaries on the tensor formalism and the definition of
some "universal" polynomials $P_m$ $(m=1,2,3,...)$
in terms of generating relations.
The rest of the paper is devoted to the proofs. In Sect.\ref{more} we
write a Leibnitz rule for the derivatives of tensor
fields on $\reali^d$, and present a generalized
Fa\`a di Bruno formula for the derivatives of any
order of $G(f, \x)$; both constructions
rely on symmetrised tensor products, which are
conveniently discussed.
In Sect.\ref{poly} we prove the main properties of the
polynomials $P_m$.
In Sect.\ref{ineq} we discuss some functional spaces of tensor fields
and the Hausdorff-Young inequality for the Fourier transform of tensors;
also, we
propose estimates for the constants in a Gagliardo type and an Adams-Frazier type
inequality. In Sect.\ref{last} we use these tools to
derive the expression of $\Upsilon_{n a d}$ presented in Sect.\ref{desc}.
Finally, in the Appendix we prove some technical statements on
tensor norms and the symmetric tensor calculus, which are necessary
to derive the main results of the paper.
\section{Description of the results.}
\label{desc}
Throughout the paper $\naturali := \{0,1,2,...\}$ and $\naturali_0 :=
\{1,2,3,...\}$; $d \in \naturali_0$ is a fixed space dimension.
\vskip 0.2cm \noindent
\textbf{Tensor notations and Sobolev spaces.} For our purposes, it is convenient
to introduce some tensor representations
for expressions with indices, and some related spaces of tensor-valued functions.
This setting will simplify the exposition of the results and the
subsequent proofs. \parn
For $m \in \naturali_0$, we indicate
with $\otimes^m \complessi^d$ the tensor product of $m$ copies
of $\complessi^d$; this will be identified with the
complex vector space of families
\beq T = (T_{\lambda})_{\lambda = (\lambda_1, ..., \lambda_m)
\in \{1, ..., d\}^m}~, \feq
where $T_{\lambda} \in \complessi$ for all $\lambda$; any such $T$
will be called a complex tensor of order $m$. The norm of
$T$ is
\beq | T | := \sqrt{ \sum_{\lambda} | T_{\lambda} |^2 }~.
\label{nt} \feq
In the sequel, we always intend $\otimes^0 \complessi^d := \complessi$
(as customary for tensor products).
Let $p \in [1, +\infty]$; then,
$L^p(\reali^d, \otimes^m \complessi^d)$ is the space of measurable
functions (modulo equality almost everywhere)
\beq T : \reali^d \vain \otimes^m \complessi^d~, \qquad x \mapsto T(x) =
(T_{\lambda}(x)) \feq
such that the function
\beq | T | : \reali^d \mapsto [0,+\infty)~, \qquad x \mapsto | T(x) | =
\sqrt{ \sum_{\lambda} | T_{\lambda}(x) |^2 } \label{modf} \feq
is in $L^p(\reali^d, \reali)$.
Obviously enough, this space is identified
with the space of families $T = (T_{\lambda})$ where $T_{\lambda}
\in L^{p}(\reali^d, \complessi)$ for each $\lambda \in \{1,...,d\}^m$.
The $L^p$ norm of a function $T$ as above is
\beq \| T \|_{L^p} := \|~ | T |~ \|_{L^p}~. \label{nlpf} \feq
We use the standard notations $\DD(\reali^d, \complessi)$,
$\DD'(\reali^d, \complessi)$ for the space
of smooth, compactly supported functions $\reali^d \mapsto \complessi$
and the space of distributions. More generally,
$\DD(\reali^d, \otimes^m \complessi^d)$ is the space of
smooth, compactly supported functions
$\Phi : \reali^d \vain \otimes^m \complessi^d$, which are identified
with the families $\Phi = (\Phi_{\lambda})_{\lambda \in \{1,...,d\}^m}$
with each component in $\DD(\reali^d, \complessi)$. We define
$\DD'(\reali^d, \otimes^m \complessi^d)$ to be the space of
families $T = (T_{\lambda})_{\lambda \in \{1,...,d\}^m}$
with each component in $\DD'(\reali^d, \complessi)$. \parn
Let $f \in C^m(\reali^d, \complessi)$, or
$f \in \DD'(\reali^d, \complessi)$. We put
\beq \nabla^m f := (\partial_{\lambda} f)_{\lambda =
(\lambda_1, ..., \lambda_m) \in \{1, ..., d\}^m} \label{weput} \feq
where $\partial_{\lambda} := \partial_{\lambda_1} ... \partial_{\lambda_m}$
($\partial_{\lambda_i}$ the partial derivative in the $\lambda_i$-th direction)
and the derivatives are intended in the ordinary or in the
distributional sense; we also intend $\nabla^0 f := f$. \parn
In the above tensor style, the $L^2$ type
space of complex Bessel potentials of integer order $n \in \naturali$
can be defined setting
\beq H^n(\reali^d, \complessi) := \{ f \in \DD'(\reali^d, \complessi)~|~
\nabla^m f \in L^2(\reali^d, \otimes^m \complessi^d)~\forall m
\in \{0, ..., n\}~\}~; \label{hnab} \feq
this space carries the (Hilbertian) norm
\beq \| f \|_n := \sqrt{\sum_{m=0}^n \left( \barray{c} n \\ m \farray \right)
\| \nabla^m f \|^{2}_{L^2}}~. \label{nonab} \feq
Definitions \rref{hnab}
and \rref{nonab} will be our standards throughout the paper.
In the integer case we consider, they coincide exactly with the conventional
definition of $H^n$ \cite{Smi} \cite{Maz} in terms
of the operator $\sqrt{1 - \Delta}^{~n}$, a power of $1$ minus the
Laplacian $\Delta$ constructed via the Fourier tranform of tempered
distributions. In particular, the norm $\| f \|_n$ in Eq.\rref{nonab}
equals $\| \sqrt{1 - \Delta}^{~n} f \|_{L^2}$, which is the
standard norm of our previous works \cite{mp1} \cite{mp2}
({\footnote{The representation of this norm via multiindices
given in \cite{mp2} for the integer case is
just a rephrasement of Eq.\rref{nonab}, less efficient for
our present purposes.}}). \parn
For $m \in \naturali$,
the tensor power $\otimes^m \reali^d$,
the spaces $L^p(\reali^d, \otimes^m \reali^d)$ and their norms are
constructed in an obvious way, with $\reali$ replacing $\complessi$ everywhere.
Similarly, we can introduce $\nabla^m f$ when $f$ is a real-valued function
or distribution on $\reali^d$; this allows to define
$H^n(\reali^d, \reali)$ and its norm writing the
analogues of Eq.s\rref{hnab} \rref{nonab}. \parn
It should be noted that, due to the permutability of partial derivatives,
the derivatives $\nabla^m f$ of a real or complex function $f$ take values in
the space of \textsl{symmetric} tensors, on which we will return later from
a more systematic viewpoint.
\vskip 0.2cm \noindent
\textbf{Useful constants and notations.} For any integer $a > d/2$, we put
\beq S_{a d} := {1 \over (4 \pi)^{d/4}} \sqrt{
{\Gamma(a - d/2) \over \Gamma(a)} }~, \label{sad} \feq
$\Gamma$ denoting the factorial function. This is the sharp
constant in the imbedding inequality of $H^{a}(\reali^d, \reali
~\mbox{or}~\complessi)$
into $L^{\infty}(\reali^d, \reali~\mbox{or}~\complessi)$; thus
\beq \|~ \|_{L^{\infty}} \leq S_{a d} \|~ \|_{a} \label{imbed} \feq
and $S_{a d}$ is the minimum constant for which this occurs \cite{mp1}.
For $a$ as before and $r \in (0, +\infty]$, we put
\beq H^{a}_r(\reali^d, \reali~\mbox{or}~\complessi) :=
\{ f \in H^{a}(\reali^d, \reali~\mbox{or}~\complessi)~|~S_{a d}
\| f \|_a < r \}~ \label{sar} \feq
(intending this to be the whole $H^a$ space, if $r = +\infty$);
then
\beq f \in H^{a}_r(\reali^d, \reali~\mbox{or}~\complessi)
\Longrightarrow \| f \|_{L^{\infty}} < r \label{long} \feq
(all this can be extended to the case of a noninteger $a > d/2$, but
here we are not considering fractional Sobolev spaces). \parn
Let $E : [0,+\infty) \vain [1,+\infty)$ be defined by
\beq E(s) := s^s \qquad \mbox{for $s \in (0, + \infty)$}~,
\qquad E(0) := \lim_{s \vain 0^{+}} E(s) = 1~. \label{defe} \feq
Using this function, we define
\beq U_{m j d} :=
\left( {E(1/2 - (j-1)/(2 m) ) \over E(1/2 + (j-1)/(2 m))} \right)^{d/2}~
\left( {E(1/(2 m)) \over E(1 - 1/(2 m))} \right)^{(j-1) d/2}~ \label{defu}
\feq
for all integers $m, j$ with $1 \leq j \leq m$. It should be noted that
\beq U_{m j d} \leq 1 \feq
for $m, j$ as above; in fact, $U_{m j d}$ is a product of factors
of the form $E(s)/E(1-s)$, with $s \in [0,1/2]$, raised to power $d/2$;
each one of these factors is $\leq 1$. We will show how the
constant $U_{m j d}$ arises in an Adams-Frazier type inequality for products
of derivatives (see Prop.\ref{adam}).
\vskip 0.2cm \noindent
\textbf{Defining the "universal" polynomials}
$\mbox{\boldmath $P_m$.}$ As anticipated in the
Introduction, these polynomials will be employed to express the main
results on Nemytskij operators; we will give some equivalent descriptions
for them via the following Lemma, to be proved in Sect.\ref{poly}. \parn
\begin{prop}
\label{lemma1}
\textbf{Lemma.} Let $m \in \naturali_0$. Then, there is a unique
real polynomial $P_m((\nu_{j \ell}),\rho)$, in a set of variables
\beq \nu_{j \ell}~,~~\rho
 \qquad \qquad(~1 \leq j \leq m~,~~0 \leq \ell \leq m - j~)~,\label{delta} \feq
such that for each $C^m$ function
$F : \reali^2 \mapsto \reali$,~$(u, \xi) \mapsto F(u, \xi)$
it is
\beq {d^m \over d \xi^m}~ F(e^\xi, \xi) = P_m
\left(\nu_{j \ell} =
{\partial^{j+\ell} F \over \partial u^j \partial \xi^\ell}(e^\xi, \xi)~,~
\rho = e^\xi \right) + {\partial^{m} F \over \partial \xi^m}(e^\xi, \xi)~.
\label{dmxm} \feq
This polynomial has the form
\beq P_m((\nu_{j \ell}), \rho) =
\sum_{\stackrel{ 1 \leq j \leq m}
{\rule{0mm}{2mm}\scrscr{0 \leq \ell \leq m-j}}} P_{m j \ell}~
\nu_{j \ell}~ \rho^j~ \label{form} \feq
where the $P_{m j \ell}$ are positive, integer coefficients; these have
the equivalent characterizations
\beq P_{m j \ell} = \left. {d^m \over d \xi^m} \right|_{\xi=0}
{(e^\xi - 1)^j \xi^\ell \over j! \ell!} ~, \label{equi1} \feq
\beq P_{m j \ell} = {1 \over j!} \left( \barray{c} m \\ \ell
\farray \right)~\sum_{s=0}^j (-1)^{j - s} \left( \barray{c} j \\ s
\farray \right) s^{m - \ell}~, \label{equi2} \feq
and can be described as the solutions of the recursive equations
$$ P_{1 1 0} = 1~, \qquad P_{1 0 0} = 0~; $$
\beq P_{m +1, j \ell} = P_{m j, \ell-1} + P_{m, j-1, \ell} + j~ P_{m j \ell} +
\delta_{j 1} \delta_{\ell m} \label{rec} \feq
$$ \mbox{for $m \geq 1$,~ $1 \leq j \leq m+1$,~ $0 \leq \ell \leq m + 1 - j$}
$$
(intending $P_{m j, -1} := P_{m, m+1, \ell} := 0$; $\delta_{i j}$ is
the Kronecker symbol).
\fine
\end{prop}
The defining relation \rref{dmxm} for the polynomials $P_m$ point
out their link with the derivation of composite functions, a topic
which has obvious connections with the analysis of Nemytskij operators;
Eq.\rref{dmxm} also allows very simple derivations of some identities
for these polynomials (see, e.g., Eq.\rref{fmth}). \parn
The automatic computation of the polynomials $P_m$ can be easily
performed, implementing Eq.s\rref{dmxm} or \rref{equi1} on any
package for symbolic manipulation, or the purely numerical
relations \rref{equi2} or \rref{rec} on any system capable to deal with integer
numbers. The polynomials $P_m$ for $1 \leq m \leq 6$
are reported in a separate Table. \parn
There is a strict relation between the coefficients
$P_{m j \ell}$
and the Stirling numbers. Let us recall that the
Stirling number of the second kind $\Stir^{j}_{m}$ is defined to be
the number of partitions of a set of $m$ elements into $j$ nonempty
subsets; comparing some standard results on these numbers \cite{Com}
with anyone of Eq.s (\ref{equi1}-\ref{rec}), it is found that
\beq P_{m j \ell} = \left( \barray{c} m \\ \ell
\farray \right) \Stir^{j}_{m - \ell} \qquad (1 \leq j \leq m,~
0 \leq \ell \leq m - j)~. \label{eqstir} \feq
(In particular, $P_{m j 0} = \Stir^{j}_m$; these are the only coefficients
involved in the derivative \rref{dmxm} in the case of a
$\xi$-independent function $F(u)$,
which is considered in \cite{Com}.
The $P_{m j 0}$ will also be the only coefficients
appearing in our estimates on
Nemystkij operators $f \mapsto G(f)$ with an $x$-independent $G$).
\parn
\begin{table}[t]
\hrule
\vskip 0.3cm \noindent
\textbf{Table of low order polynomials} $\mbox{\boldmath $P_m$.}$
\vskip 0.3cm \noindent
$ \dd
P_1((\nu_{j \ell}),   \rho) = \nu_{1  0} \rho$
\vskip 0.1cm \noindent
$ \dd P_2((\nu_{j \ell}),   \rho) =
\nu_{2  0} \rho^2 + (\nu_{1  0} + 2 \nu_{1  1}) \rho$
\vskip 0.1cm \noindent
$ \dd P_3((\nu_{j \ell}),   \rho) =
\nu_{3  0} \rho^3 +
(3 \nu_{2  0} + 3 \nu_{2  1}) \rho^2 +
(\nu_{1  0} + 3 \nu_{1  1} + 3 \nu_{1  2}) \rho$
\vskip 0.1cm \noindent
$ \dd P_4((\nu_{j \ell}),   \rho) =
\nu_{4  0} \rho^4 +
(6 \nu_{3  0} + 4 \nu_{3  1}) \rho^3
+ ( 7 \nu_{2  0} + 12 \nu_{2  1} + 6 \nu_{2  2}) \rho^2 $
$ \dd
+ (\nu_{1  0} + 4 \nu_{1  1} + 6 \nu_{1  2} + 4 \nu_{1  3}) \rho$
\vskip 0.1cm \noindent
$ \dd
P_5((\nu_{j \ell}), \rho) = \nu_{5 0} \rho^5 +
(10 \nu_{4 0} + 5 \nu_{4 1}) \rho^4 +
(25 \nu_{3 0} + 30 \nu_{3 1} + 10 \nu_{3 2}) \rho^3 + $ \parn
$ \dd
(15 \nu_{2 0} + 35 \nu_{2 1} + 30 \nu_{2 2} + 10 \nu_{2 3}) \rho^2 +
(\nu_{1 0} + 5 \nu_{1 1} + 10 \nu_{1 2} + 10 \nu_{1 3}
+ 5 \nu_{1 4}) \rho $
\vskip 0.1cm \noindent
$ \dd
P_6((\nu_{j \ell}), \rho) = \nu_{6 0} \rho^6 +
(15 \nu_{5 0} + 6 \nu_{5 1}) \rho^5 +
(65 \nu_{4 0} + 60 \nu_{4 1} + 15 \nu_{4 2}) \rho^4 + $ \parn
$ \dd
+ (90 \nu_{3 0} + 150 \nu_{3 1} + 90 \nu_{3 2} + 20 \nu_{3 3}) \rho^3
+ (31 \nu_{2 0} + 90 \nu_{2 1} + 105 \nu_{2 2}
+ 60 \nu_{2 3} + 15 \nu_{2 4}) \rho^2 + $ \parn
$ \dd +
(\nu_{1 0} + 6 \nu_{1 1} + 15 \nu_{1 2} + 20 \nu_{1 3} + 15 \nu_{1 4} +
 6 \nu_{1 5}) \rho $
\vskip 0.3cm \noindent
\hrule
\end{table}
\vskip 0.2cm \noindent
\textbf{Nemytskij operators in
the real case}. Let us consider a function
\beq G : (- r, r) \times \reali^d \vain \reali~, \qquad
(u, x) = (u, x_1, ..., x_d) \mapsto G(u, x) \feq
with $r \in (0, +\infty]$. Concerning partial derivatives in the above
variables, we use the notations
\beq \partial^{j} := {\partial^j \over \partial u^j}
\qquad (j \in \naturali)~,
\feq
\beq
\partial_{\lambda} := \partial_{\lambda_1} ... \partial_{\lambda_\ell} :=
{\partial \over \partial x_{\lambda_1}} ....
{\partial \over \partial x_{\lambda_\ell}} \qquad (\ell \in \naturali,
\lambda = (\lambda_1, ..., \lambda_\ell) \in \{1, ..., d\}^\ell)~.
\label{intend} \feq
For $j, \ell \in \naturali$ and $G$ of class $C^{j + \ell}$, the family
of partial derivatives
\beq \partial^{j} \nabla^\ell G :=
\left(\partial^j \partial_{\lambda} G \right)_{\lambda \in \{1,...,d\}^\ell}
\label{nhnk} \feq
can be seen as a continuous function
$\partial^j \nabla^\ell G : (-r, r)
\times \reali^d \vain \otimes^\ell \reali^d$.
\begin{prop}
\label{defuno}
\textbf{Definition.} Given the map $G : (-r,r) \times \reali^d
\vain \reali$, we put
\beq \flat_m(G,\rho) :=
\sup_{\barray{c} \scr{u \in [-\rho,\rho] \setminus \{0 \} } \\
\scr{x \in \reali^d} \farray }
\left| {\nabla^m G(u, x) - \nabla^m G(0, x) \over u} \right|
\qquad (~\rho \in (0, r)~)~,
\label{sim} \feq
\beq \sigma_{j \ell}(G,\rho) :=
\sup_{\barray{c} \scr{u \in [-\rho,\rho]} \\
\scr{x \in \reali^d} \farray }
\left| \partial^j \nabla^\ell G(u, x) \right|
\qquad (~\rho \in [0, r)~)
\label{sihk} \feq
if $m, j, \ell$ are integers and $G$ is $C^{m}$ or $C^{j+\ell}$, respectively
(the indicated suprema being possibly infinite). Also, we put
\beq \beta_{m d}(G,\rho) :=
P_m \left(\nu_{j \ell} = (1 - {\ell \over m}) U_{m - \ell, j d}
{\sigma_{j \ell}(G,\rho) \over \rho},~\rho \right)
\qquad (~\rho \in (0, r)~)~,
\label{betam} \feq
\beq b_{m d}(G,\rho) := P_m \left(\nu_{j \ell} =
{\ell \over m} U_{m - \ell, j  d} {\sigma_{j \ell}(G,\rho) \over \rho},~
\rho \right)
\qquad (~\rho \in (0, r)~)~,
\label{bm} \feq
if $m \geq 1$ is an integer and $G$ is $C^m$ (the r.h.s. of
the last two equations being possibly
$+\infty$, when infinite suprema $\sigma_{j \ell}(G, \rho)$
are involved; in any case,
it must be understood that $\nu_{j 0} =0$ in Eq.\rref{bm}).
Here, $P_m$ is the polynomial defined by Lemma
\ref{lemma1} and the constants $U_{m - \ell,j d}$ are
defined according to Eq.\rref{defu}.
\fine
\end{prop}
Of course $\flat_m(G, \rho) \leq \flat_m(G, \rho_1)$,
$\sigma_{j \ell}(G,\rho) \leq \sigma_{j \ell}(G,\rho_1)$ for
$\rho \leq \rho_1$; thus, if anyone of these suprema is finite
for some radius $\rho_1$, it is also finite for all radii $\rho \leq \rho_1$.
For uniformity of language, it is convenient
to define $\flat_m(G,\rho)$, $\beta_{m d}(G, \rho)$ and
$b_{m d}(G, \rho)$ at $\rho = 0$ also. \parn
If $\flat_m(G, \rho)$ is finite for $\rho$ in a right neighborhood
of zero, we will make the natural choice
$\flat_m(G,0) :=$ $\lim_{\rho \vain 0^+} \flat_m(G, \rho)$.
By monotonicity this limit exists; it is finite and $\geq 0$.
\parn
We will also intend $\beta_{m d}(G,0) := \lim_{\rho \vain 0^+}
\beta_{m d}(G,\rho)$, $b_{m d}(G,0) := \lim_{\rho \vain 0^+}
b_{m d}(G,\rho)$, if the r.h.s. are finite for $\rho$ close
to zero. Again, the limits exist finite
and $\geq 0$: in fact, in the $\nu_{j \ell}$ arguments of
$P_m$ the terms $\sigma_{j \ell}(G, \rho)$
admit finite limits by monotonicity, and the factors $1/\rho$ are
cancelled out by the powers $\rho^j$ ($j \geq 1$), see Eq.\rref{form}.
\vskip 0.2cm
\begin{prop}
\label{defunoa}
\textbf{Definition.} Let $n \in \naturali$.
A map $G$ as above is said to have the $\Pi^n$ property if: \parn
i) $G$ is $C^n$; \parn
ii) $\flat_m(G, \rho) < +\infty$ for $0 \leq m \leq n$ and
all $\rho \in (0, r)$;
$\sigma_{j \ell}(G, \rho) < + \infty$ for $1 \leq j \leq n$,
$0 \leq \ell \leq n -j$ and all $\rho \in [0,r)$; \parn
iii) the functions $\flat_m(G, .)$,
$\sigma_{j \ell}(G, .) : [0, r) \vain [0,+\infty)$
are continuous for all $m, j, \ell$ as in ii).
\fine
\end{prop}
Of course, iii) implies continuity of the functions $b_{m d}(G, .),
\beta_{m d}(G, .) : [0, r) \vain [0, +\infty)$ for $1 \leq m \leq n$.
Let us recall the definitions \rref{sad}
of the constant $S_{a d}$ and \rref{sar} of the ball
$H^a_{r}(\reali^d, \reali)$, where $a > d/2$;
for $G : (-r, r) \times \reali^d \vain \reali$
and $f \in H^a_{r}(\reali^d, \reali)$, $G(f, \x)$ is defined
on the grounds of Eq.\rref{long}. \parn
With the above terminology, the main result
on real Nemystkij operators is the following.
\begin{prop}
\label{nemreale}
\textbf{Proposition.} Let $n, a \in \naturali$ and $a > d/2$. Consider: \parn
i) a function
$G : (-r, r) \times \reali^d \vain \reali$ with the $\Pi^n$ property,
such that $G(0,\x) : x \mapsto G(0, x)$ is in $H^{n}(\reali^d, \reali)$; \parn
ii) a function $f \in H^{n}(\reali^d, \reali) \cap
H^a_{r}(\reali^d, \reali)$. \parn
Then $G(f,\x) \in H^{n}(\reali^d, \reali)$ and
\beq \| G(f,\x) - G(0,\x) \|_n \leq
\gamma_{n d}(G, S_{a d} \| f \|_a) \| f \|_n +
c_{n d}(G, S_{a d} \| f \|_a) \| f \|_{L^2}
\label{disug} \feq
where, for any $\rho \in [0,r)$,
\beq \gamma_{n d}(G, \rho) :=
\max \left\{ \flat_{0}(G, \rho),~\beta_{m d}(G, \rho)~
(m=1,...,n) \right\} ,  \label{gammand} \feq
\beq c_{n d}(G, \rho) :=
\sqrt{\sum_{m=1}^n \left( \barray{c} n \\ m \farray \right)
\left( b_{m d}(G, \rho) + \flat_m(G, \rho) \right)^2} \label{cnd} \feq
(intending these to be $\flat_0(G,\rho)$ and $0$, respectively,
for $n=0$). Eq.\rref{disug} implies the weaker bound
\beq \| G(f,\x) - G(0,\x) \|_n \leq
\Big(~ \gamma_{n d}(G, S_{a d} \| f \|_a) +
c_{n d}(G, S_{a d} \| f \|_a)~ \Big)~\| f \|_n~. \label{weak} \feq
\fine
\end{prop}
\textbf{Remarks. i)} Of course, the coefficient of $\| f \|_n$
in the r.h.s. of Eq.\rref{weak} is the term
denoted with $\Upsilon_{n a d}(G, \| f \|_a)$
in the Introduction. \parn
\textbf{ii)} We recall that $\beta_{m d}$ and $b_{m d}$ depend on the constants
$ U_{m - \ell, j  d}~$, all of them $\leq 1$; weaker but
simpler estimates can be obtained replacing systematically
those constants with one. \parn
\textbf{iii)} The previous results become simpler
if $G$ is $x$-independent: $G(u, x) = G(u)$; this case is frequently
considered in the investigation of Nemytskij operators and will be discussed
in more detail at the end of this Section, also giving two examples.
A main feature of this case is that $c_{n d}(G, \rho) = 0$ for all
$n$, $\rho$. \fine
\vskip 0.2cm\noindent
\textbf{Nemytskij operators in the complex case}.
From now on $\B(0,r)$
denotes the open ball in $\complessi$
of center $0$ and radius $r$, and $\Bc(0,r)$ its closure.
We consider a function
\beq G : \B(0,r) \times \reali^d \vain \complessi~, \qquad
(z, x) = (z, x_1, ..., x_d) \mapsto G(z, x) \feq
with $r \in (0, +\infty]$. For partial derivatives with respect to
the complex variable, we use the notations
\beq \partial^{h} \ppartial^{k} := {\partial^h \over \partial z^{h}}
{\partial^{k} \over \partial \z^{k}}
\qquad (h,k \in \naturali)~;
\feq
we write $\partial_{\lambda}$, as in Eq.\rref{intend}, for the derivatives
w.r.t. the $x$ variables. We say that $G$ is $C^n$ if
all partial derivatives of $G$ of order
$\leq n$ w.r.t. $z, \overline{z}$ and $(x_1, ..., x_d)$
exist and are continuous, everywhere on its domain.  We put
\beq \partial^{h} \ppartial^{k} \nabla^\ell G :=
\left(\partial^h \ppartial^{k}
\partial_{\lambda} G \right)_{\lambda \in \{1,...,d\}^\ell}
\label{nhnkc} \feq
for $h, k, \ell \in \naturali$ and $G$ of class $C^{h+k+\ell}$;
this can be regarded as a function
$\B(0,r) \times \reali^d \vain \otimes^\ell \complessi^d$.
\begin{prop}
\label{defunoc}
\textbf{Definition.} Given the function $G : \B(0, r) \times \reali^d
\vain \complessi$, we put
\beq \flat_m(G,\rho) :=
\sup_{\barray{c} \scr{z \in \Bc(0,\rho) \setminus \{0 \} } \\
\scr{x \in \reali^d} \farray }
\left| {\nabla^m G(z, x) - \nabla^m G(0, x) \over z} \right|
\qquad (~\rho \in (0, r)~)~,
\label{simc} \feq
\beq \sigma_{h k \ell}(G,\rho) :=
\sup_{\barray{c} \scr{z \in \Bc(0,\rho)} \\
\scr{x \in \reali^d} \farray }
\left| \partial^h \ppartial^{k} \nabla^\ell G(z, x) \right|
\qquad (~\rho \in [0, r)~)
\label{sihkc}~, \feq
\beq \sigma_{j \ell}(G,\rho) :=
\sum_{h=0}^j \left( \barray{c} j \\ h \farray \right)~
\sigma_{h, j-h, \ell}(G, \rho)
\qquad (~\rho \in [0, r)~) \label{sijlc} \feq
for $m, h, k, j, \ell$ integers and $G$ of class $C^m$,
$C^{h+k+\ell}$ and $C^{j+\ell}$, respectively (these quantities
being possibly infinite). \parn
For $m \geq 1$ integer and $\rho \in (0, r)$,
the functions $\beta_{m d}(G, \rho)$, $b_{m d}(G, \rho)$
will be defined as in Def.\ref{defuno}, using again the
polynomials $P_m$ of Lemma \ref{lemma1}.
\fine
\end{prop}
As in the real case, the indicated suprema are monotonically increasing with
$\rho$, and their finiteness at $\rho_1$ ensures finiteness for all
$\rho \leq \rho_1$. If $\flat_m(G, \rho) < +\infty$, at least
for $\rho$ in a right neighborhood of zero, we put
$\flat_m(G, 0) :=$ $\lim_{\rho \vain 0^+} \flat_m(G, \rho)$.
We also intend $\beta_{m d}(G,0) := \lim_{\rho \vain 0^+}
\beta_{m d}(G,\rho)$, $b_{m d}(G,0) := \lim_{\rho \vain 0^+}
b_{m d}(G,\rho)$, if these quantities are finite for $\rho$ close
to zero.
\parn
\begin{prop}
\label{defunoca}
\textbf{Definition.} Let $n \in \naturali$. A map $G$ as above is
said to have the $\Pi^n$ property if: \parn
i) $G$ is $C^n$; \parn
ii) $\flat_m(G, \rho) < +\infty$ for $0 \leq m \leq n$ and
all $\rho \in (0, r)$;
$\sigma_{h k \ell}(G, \rho) < + \infty$ for $h, k \geq 0$,
$1 \leq h+k \leq n$, $0 \leq \ell \leq n -h-k$
and all $\rho \in [0,r)$; \parn
iii) the functions $\flat_m(G, .)$,
$\sigma_{h k \ell}(G, .) : \B(0,r) \vain [0,+\infty)$
are continuous for all $m, h, k, \ell$ as in ii).
\fine
\end{prop}
The main result on complex Nemytskij operators is the following.
\begin{prop}
\label{nemc}
\textbf{Proposition.} Let $n, a \in \naturali$ and
$a > d/2$. Consider functions $G : \B(0,r) \times \reali^d \vain \complessi$,
$f : \reali^d \vain \complessi$ with the
$\complessi$-analogues of properties i) ii) in
Prop.\ref{nemreale}.
Then $G(f,\x) \in H^{n}(\reali^d, \complessi)$ and Eq.s
(\ref{disug}-\ref{weak}) hold again, with
the definitions given therein for $\gamma_{n d}(G,\rho)$ and
$c_{n d}(G,\rho)$.
\fine
\end{prop}
\vskip 0.2cm \noindent
\textbf{The case of an $\mbox{\boldmath $x$}$-independent
$\mbox{\boldmath $G$}$.} Let us consider a function
\beq G : (- r, r) \vain \reali,~
u \mapsto G(u) \qquad \mbox{or} \qquad
G : \B(0,r) \vain \complessi,~ z \mapsto G(z) \feq
(derivable in the ordinary sense, as many times as required in the sequel);
we regard $G$ as an $x$-independent function of $(u,x)$ or $(z, x)$.
Some of the previous formulae become
simpler due to the vanishing of all derivatives $\nabla^m G$ $(m \geq 1)$
with respect to the $x$ variables. More precisely,
\beq \flat_m(G, \rho) = 0
\quad \mbox{for $m \geq 1$}~, \qquad
\sigma_{j \ell}(G, \rho) = 0 \quad \mbox{for $\ell \geq 1$}~, \feq
also implying
\beq b_{m d}(G, \rho) = 0 \qquad \mbox{for $m \geq 1$}~, \feq
\beq \beta_{m d}(G,\rho) =
P_m \left(\nu_{j 0} = U_{m j d}
{\sigma_{j 0}(G,\rho) \over \rho};
\nu_{j \ell} = 0~\mbox{for $\ell \neq 0$}~,
~\rho \right)~; \label{must} \feq
\beq c_{n d}(G, \rho) = 0 \qquad
\mbox{for all $n \geq 0$}~. \feq
$G$ has the $\Pi^n$ property for $n \geq 1$ if and only if it is $C^n$;
the $\Pi^0$ property is possessed by $G$ if and only if $G$ is $C^0$
and $\flat_0(G, \rho)$ is finite and continuous in $\rho$. \parn
Both in the real and in the complex case,
we can study the Nemytskij operator sending
$f \in H^n(\reali^d) \cap H^{a}_r(\reali^d)$ into $G(f)$,
with the aid of Prop.s \ref{nemreale} or \ref{nemc}. Of course,
the constant function $x \in \reali^d \mapsto G(0)$ is in
$H^n(\reali^d)$ for any $n$ if and only if
\beq G(0) = 0~. \feq
If this happens and
$G$ has the $\Pi^n$ property, Eq.\rref{disug} holds, taking in this case
the form
\beq \| G(f)  \|_n \leq
\gamma_{n d}(G, S_{a d} \| f \|_a) \| f \|_n ~; \label{disugnox} \feq
as in Eq.\rref{gammand}, it is
$\gamma_{n d}(G, \rho) :=
\max \left\{ \flat_{0}(G, \rho),~\beta_{m d}(G, \rho)~
(m=1,...,n) \right\}$. We shall exemplify this scheme in two cases:
in the first $G$ is a monomial in $z$ and $\overline{z}$,
in the second it is the hyperbolic sine.
\vskip 0.2cm \noindent
\textbf{Examples.}
\textbf{i)} Let us put
\beq G : \complessi \vain \complessi~, \qquad z \mapsto G(z) := z^{\H}
\overline{z}^{\K}~, \feq
where $\H, \K \in \naturali$ and $\H +\K \neq 0$.
Of course $\partial^h \ppartial^{k} G(z) =
\H! \K! /(\H - h)! (\K- k)!$ $z^{\H-h} \overline{z}^{\K -k}$ for
$0 \leq h \leq \H$
and $0 \leq k \leq \K$; all the other derivatives of $G$ vanish. The nonzero
functions in Eq.s\rref{simc} \rref{sihkc} are
\beq \flat_0(G, \rho) = \rho^{\H +\K -1}~, \feq
\beq \sigma_{h k 0}(G, \rho) = {\H! \K! \over (\H - h)! (\K - k)!}~
\rho^{\H +\K -h- k} \qquad (0 \leq h \leq \H; 0 \leq k \leq \K)~; \feq
from here, one also infers that the only nonzero
functions in \rref{sijlc} are the $\sigma_{j 0}$ for
$1 \leq j \leq \H + \K$, given by
\beq \sigma_{j 0}(G, \rho) = \left( \sum_{h=0}^j
\left( \barray{c} j \\ h \farray \right)
{\H! \K! \over (\H - h)! (\K - j + h)!} \right)~\rho^{\H +\K -j} = \feq
$$  = j! \left( \sum_{h=0}^j \left( \barray{c} \H \\ h \farray \right)~
\left( \barray{c} \K \\ j - h \farray \right) \right)~\rho^{\H +\K -j} =
j! \left( \barray{c} \H + \K \\ j \farray \right)~\rho^{\H + \K - j} =$$
$$ ={(\H + \K)! \over (\H + \K - j)!}~\rho^{\H + \K - j}. $$
Thus, we have
$$ \beta_{m d}(G,\rho) = P_m \Big(\nu_{j 0} =
{(\H + \K)! \over (\H + \K - j)!}~U_{m j d}~ \rho^{\H + \K - j -1}
~\mbox{for $1 \leq j \leq \H + \K$}~;~$$
$$ \nu_{j \ell} = 0~\mbox{otherwise}~,~ \rho \Big)~. $$
On the other hand,
each coefficient $\nu_{j 0}$
in the polynomial $P_m$ multiplies $\rho^{j}$
(see Eq.\rref{form}). From here,
one easily infers
\beq \beta_{m d}(G, \rho) = B^{\H + \K}_{m d}
\rho^{\H + \K -1} \qquad (m \geq 1) \feq
where the $B^{\J}_{m d}$ are positive coefficients defined
for any integer $\J \geq 1$, setting
\beq B^{\J}_{m d} :=
P_m\Big(\nu_{j 0} =
{\J! \over (\J - j)!}~U_{m j d}~\mbox{for $1 \leq j \leq \J$}~;~
 \nu_{j \ell} = 0~\mbox{otherwise}~,~ \rho = 1\Big)~; \feq
this also implies
\beq \gamma_{n d}(G, \rho) = \Gamma^{\H + \K}_{n d} \rho^{\H +\K -1}
\label{gammahk} \feq
$$ \Gamma^{\J}_{n d} :=
\max\left(1,~ B^{\J}_{m d}~(m=1,...,n) \right) \qquad (n \in \naturali)~.
$$
The final result on the Nemytskij operator
$G(f) = f^{\H} \overline{f}^{\K}$ is the following, for all $n, a \in \naturali$
with $a > d/2$~: for each
$f \in H^n(\reali^d, \complessi) \cap H^a(\reali^d, \complessi)$ it is
$f^{\H} \overline{f}^{\K} \in H^n(\reali^d, \complessi)$ and
\beq \| f^{\H} \overline{f}^{\K} \|_{n} \leq
\Gamma^{\H + \K}_{n d} (S_{a d} \| f \|_a)^{\H+\K-1} \| f \|_n~. \label{tha} \feq
A slightly weaker, but simpler bound can be obtained recalling that
all the $U$ coefficients are $\leq 1$.
For each $m \geq 1$, this implies
$$ B^{\J}_{m d} \leq
P_m\Big(\nu_{j 0} =
{\J! \over (\J - j)!}~\mbox{for $1 \leq j \leq \J$}~;~
\nu_{j \ell} = 0~\mbox{otherwise}~,~ \rho = 1\Big) = $$
\beq
= \left. {d^m \over d \xi^m} \right|_{\xi=0} e^{\J \xi} = \J^m~;
\label{fmth} \feq
the equality at the beginning of the last line can be obtained directly
from the relation \rref{dmxm} defining $P_m$, to be applied at $\xi=0$ to
the function $F(u, \xi) := u^{\J}$. \parn
From \rref{gammahk} and \rref{fmth}
one readily infers $\Gamma^{\J}_{n d} \leq$ $\J^n$ for all $n \geq 0$,
whence
\beq \| f^{\H} \overline{f}^{\K} \|_{n} \leq
(\H + \K)^n (S_{a d} \| f \|_a)^{\H +\K -1} \| f \|_n~. \label{thb} \feq
Let us recall that, in the previous work \cite{mp2}, we
have discussed some inequalities for the norm $\| f g \|_n$ of a
product, implying
directly an estimate of the form
$\| f^{\H} \overline{f}^{\K} \|_n \leq K_{n a d} \| f \|_a \| f \|_n$
for all $\H$, $\K$ with $\H +\K=2$~
({\footnote{In \cite{mp2}, it was assumed either $0 \leq n \leq d/2
< a$ or $d/2 < a \leq n$.}}). In all the tests we have performed,
the upper bounds
of \cite{mp2} on the
constants $K_{n a d}$ and the
bounds $\Gamma^{\H + \K}_{n d}
S_{a d}$ arising from \rref{tha}, with $\H + \K = 2$,
are very close numerically.
A result of \cite{mp2} also implies the weaker bound
\rref{thb} for $\H + \K=2$. \parn
For $\H + \K > 2$, an estimation of $\| f^{\H} \overline{f}^{\K} \|_{n}$
could be performed by iteration of the bounds \cite{mp2}
on the product of two functions; however, the estimates derived in this
elementary way would be sensibly rougher than \rref{tha} or \rref{thb}.
\vskip 0.2 cm \noindent
\textbf{ii)} We consider the hyperbolic function
\beq G : \reali \vain \reali~, \qquad u \mapsto G(u) := \sinh u~. \feq
From Eq.s\rref{sim} \rref{sihk} we get
\beq \flat_0(\sinh, \rho) = {\sinh \rho \over \rho}~, \qquad
\flat_m(\sinh, \rho) = 0
\qquad (m \geq 1) \feq
\beq \sigma_{j 0}(\sinh, \rho) = \sinh \rho \quad \mbox{($j$ even)},
\qquad
\sigma_{j 0}(\sinh, \rho) = \cosh \rho \quad \mbox{($j$ odd)}~, \feq
$$ \sigma_{j \ell}(\sinh, \rho) = 0 \qquad \mbox{for $\ell \neq 0$}~. $$
From these objects and Eq.\rref{betam} one computes
\beq \beta_{m d}(\sinh, \rho) =
P_m \Big(\nu_{j 0} = U_{m j d}
{\sinh \rho \over \rho}~\mbox{for $j$ even}~; \feq
$$ \nu_{j 0} = U_{m j d}
{\cosh \rho \over \rho}~\mbox{for $j$ odd}~;
\nu_{j \ell} = 0~\mbox{for $\ell \neq 0$}~, ~\rho \Big)~ $$
and $\gamma_{n d}(\sinh, \rho) =
\max(\flat_0(\sinh, \rho),
\beta_{1 d}(\sinh, \rho), ... ,\beta_{n d}(\sinh, \rho)$). We have
\beq \| \sinh f \|_n \leq
\gamma_{n d}(\sinh , S_{a d} \| f \|_a) \| f \|_n \label{sinip} \feq
for all $n, a$ as usually and $f \in H^n(\reali^d, \reali) \cap H^a(\reali^d,
\reali)$. For instance, let us choose $n=2$;
this requires use of the polynomials $P_1$, $P_2$ in the Table,
from which one computes
$\beta_{1 d}(\sinh, \rho) = U_{1 1 d} \cosh \rho$ and
$\beta_{2 d}(\sinh, \rho) =
U_{2 2 d}~ \rho \sinh \rho + U_{2 1 d} \cosh \rho$. Eq.\rref{sinip}
holds with
\beq \gamma_{2 d}(\sinh, \rho) = \max \Big({\sinh \rho \over \rho},
 U_{1 1 d} \cosh \rho,~
U_{2 2 d}~ \rho \sinh \rho + U_{2 1 d} \cosh \rho~\Big).\feq
\section{More on tensor calculus. Iterated derivatives of Nemytskij
operators.}
\label{more}
\salto
The purpose of this Section and of the next two
is to give some tools for proving the main results of the paper.
Here we concentrate on some algebraic aspects of tensor calculus, on
the differentiation of tensors and on a general formula for the
derivatives of any order of a composite function $G(f,\x)$.
In Sect.\ref{poly} we will connect these results with the
polynomials $P_m$; in Sect.\ref{ineq} we will be concerned with functional
spaces of tensor-valued mappings, and some related inequalities. \parn
We already introduced the tensor powers $\otimes^\ell \complessi^d$. For
$\ell, m \in \naturali$, the tensor product operation
\beq \otimes : (\otimes^\ell \complessi^d) \times (\otimes^m
\complessi^d) \vain \otimes^{\ell + m} \complessi^d~,
\qquad (T, U) \mapsto T \otimes U \feq
is defined setting
\beq (T \otimes U)_{\lambda_1 ... \lambda_\ell \mu_1 ...
\mu_m} := T_{\lambda_1 ... \lambda_\ell} U_{\mu_1 ... \mu_m} \feq
for $\lambda_1, ..., \mu_m \in \{1, ..., d\}$ (for
$\ell=0$, $T \otimes U$ means the product $T U$ between the complex
number $T$ and $U$, and the same is understood for $m=0$). The product
$\otimes$ is clearly associative. \parn
Let $\ell \in \naturali_0$
and $\sigma$ a permutation of $\{1, ..., \ell\}$; then, we a have
a (linear) permutation operator
\beq \PP_{\sigma} : \otimes^{\ell} \complessi^d \vain
\otimes^{\ell} \complessi^d~, \qquad
(\PP_{\sigma} T)_{\lambda_1 ... \lambda_\ell} :=
T_{\lambda_{\sigma(1)} ... \lambda_{\sigma(\ell)}}~.
\label{perm} \feq
A tensor $T \in \otimes^{\ell} \complessi^{d}$ is
said to be symmetric if $\PP_{\sigma} T := T$ for all permutations
$\sigma$; any $T \in \otimes^{0} \complessi^d = \complessi$ is
symmetric by definition. Let $\ell \in \naturali$; we put
\beq \vee^{\ell} \complessi^d := \{ T \in \otimes^{\ell} \complessi^d~|~
\mbox{$T$ is symmetric}~\} \feq
and introduce a symmetrisation operator
\beq \SS : \otimes^{\ell} \complessi^d \vain
\vee^{\ell} \complessi^d~, \qquad
\SS := {1 \over \ell!} \sum_{\sigma \in \boma{\ell!}}
\PP_{\sigma} \label{simm} \feq
where $\boma{\ell!}$ stands for the set of all
permutations of $\{1, ..., \ell\}$
(whose cardinality is $\ell!$).
Note that $\SS$ is the identity map for $\ell=1$; we
intend $\SS$ to be the identity for $\ell=0$ also. A tensor $T$ is
symmetric if and only if $\SS T = T$. \parn
For $\ell, m \in \naturali$, the symmetrised tensor product
operation $\vee$ is defined setting
\beq \vee : (\otimes^\ell \complessi^d) \times (\otimes^m
\complessi^d) \vain \vee^{\ell + m} \complessi^d~,
\qquad (T, U) \mapsto T \vee U := \SS(T \otimes U)~. \feq
This operation is commutative and associative (see the Appendix) :
\beq U \vee T = T \vee U~,
\qquad (T \vee U) \vee V = T \vee (U \vee V) \label{assoc} \feq
for arbitrary tensors $T, U$ as above and $V
\in \otimes^p \complessi^d$. (Often, $\vee$
is considered to act on symmetric tensors only; this leads to
the so-called symmetric tensor algebra, see, e.g., \cite{Ster} \cite{Kos}).
If $T \in \otimes^\ell \complessi^d$
and $q \in \naturali$, we will write
\beq \vee^q T := T \vee T ... \vee T \quad \mbox{($q$ times)}~, \label{veeq} \feq
intending this to be $1$ if $q = 0$. \parn
We already defined the norm $|~|$ of tensors, see Eq.\rref{nt};
for all $T, U$ it is
\beq | T \otimes U | = | T | | U |~, \qquad | T \vee U | \leq | T | | U |
\label{notv} \feq
(see the Appendix).
Sometimes, we need to employ the complex conjugate $\overline{T}$ of a tensor
$T \in \otimes^\ell \complessi^d$; this has components
\beq \overline{T}_{\lambda_1 ... \lambda_\ell} :=
\overline{T_{\lambda_1 ... \lambda_\ell}}~. \label{conj} \feq
We now pass to spaces of tensor-valued functions, or tensor fields
\beq T : \reali^d \mapsto \otimes^{\ell} \complessi^d~\mbox{or}~
\vee^\ell \complessi^d, \qquad
x = (x_1, ..., x_d) \mapsto T(x)~; \label{asab} \feq
the products $\otimes, \vee$, the conjugate
${~}^{\overline{~}}$, etc. are defined pointwisely on such
functions. $C^{n}(\reali^d, \otimes^{\ell} \complessi^d~\mbox{or}~
\vee^\ell \complessi^d)$ is the space of mappings as in \rref{asab},
where all components are $C^n$ in the usual sense. \parn
Given a tensor field $T$ as above, of class $C^1$, we define the derivative
$\nabla T$ to be the tensor field of order $\ell + 1$ with components
\beq (\nabla T)_{\lambda_1 ... \lambda_{\ell+1}} :=
\partial_{\lambda_{\ell+1}} T_{\lambda_1 ... \lambda_\ell}~.
\label{nab} \feq
For $p \in \naturali$
we denote with $\nabla^p$ the $p$-th power of
the operator $\nabla$ (acting on $C^p$ tensor fields); for any
$T \in C^p(\reali^d, \otimes^\ell \complessi^d)$, the tensor field
$\nabla^p T$ $\in C(\reali^d, \otimes^{\ell+p} \complessi^d)$ has components
\beq (\nabla^p T)_{\lambda_1...\lambda_{\ell+p}} =
(\partial_{\lambda_{\ell+p}} ...
\partial_{\lambda_{\ell+1} }) T_{\lambda_1...\lambda_\ell} \feq
(of course derivatives $\partial_{\lambda_{\ell+p}}$, ...
$\partial_{\lambda_{\ell+1}}$ commute, so their ordering is
immaterial).
In the sequel, we often employ the symmetrised derivative
\beq \nabla_{\SS} T := \SS(\nabla T)~, \feq
and denote with $\nabla_{\SS}^p$ the $p$-th power of
the operator $\nabla_{\SS}$.
On any function $f \in C^p(\reali^d, \complessi)$, we have
\beq \nabla^p_{\SS} f = \nabla^p f = \mbox{the already defined tensor field
\rref{weput}}~. \feq
One can write a Leinitz rule for $\nabla$ and the
product $\otimes$, see the Appendix; we will not use
this rule directly and always refer to a symmetrized
version, also proved in the Appendix.
The "symmetrized" Leibnitz rule reads
\beq \nabla_{\SS}(T \vee U) =
(\nabla_{\SS} T) \vee U + T \vee (\nabla_{\SS} U)~  \label{lei2} \feq
for all $C^1$ tensor fields $T : \reali^d \vain \otimes^\ell
\complessi^d$ and $U : \reali^d \vain \otimes^m \complessi^d$.
\parn
What has been said up to now for complex tensors can be rephrased
in the real case; so, we can define the product $\otimes$
between tensors $T \in \otimes^\ell \reali^d$, $U \in \otimes^{m}
\reali^d$, the symmetrisation operator $\SS$, the $\vee$ product
and the derivatives $\nabla, \nabla_{\SS}$ of real
tensor fields, which satisfy again the previously mentioned
Leibnitz rules.
\vskip 0.2cm \noindent
\textbf{Parameter dependent tensor fields.} We begin
from the real case. We consider tensor fields depending on a
supplementary real variable $u$, i.e., mappings
\beq T : (-r,r) \times \reali^d \vain \otimes^{\ell} \reali^d~, \qquad
(u, x) \mapsto T(u, x)~; \label{asin} \feq
we define all previous operations $\otimes$, $\SS$, $\vee$, $\nabla$,
$\nabla_{\SS}$ treating $u$ as a parameter, e.g.,
$(T \otimes U)(u,x) := T(u,x) \otimes U(u,x)$,
$(\nabla T)_{\lambda_1 ... \lambda_{\ell+1}}(u,x) :=$
$\partial_{\lambda_{\ell+1}} T_{\lambda_1 ... \lambda_\ell}(u,x)$
with $\partial_{\lambda_i} := \partial/\partial x_{\lambda_i}$.
As in Sect.\ref{desc}, we denote with $\partial$ the derivative
$\partial/\partial u$. For any $u$-dependent tensor field $T$
of order $\ell$ as in \rref{asin} and class $C^1$,
$\partial T$ is the $u$-dependent tensor field, again of order $\ell$,
with components
\beq (\partial T)_{\lambda_1 ... \lambda_{\ell}} :=
\partial (T_{\lambda_1 ... \lambda_{\ell}})~. \feq
Of course $\partial$ commutes with  $\nabla$ and $\nabla_{\SS}$.
We can iterate the operations $\partial, \nabla, \nabla_{\SS}$; so,
for $T$ of order $\ell$ and class $C^{j+m}$, the tensor fields
$\partial^j \nabla^m T$, $\partial^j \nabla_{\SS}^m T$
(both of order $\ell+m$) are defined. \parn
The case of complex tensor fields, depending a complex parameter
$z$, is worked out similarly. In this case we consider mappings
\beq T : \B(0,r)
\times \reali^d \vain \otimes^{\ell} \complessi^d~, \qquad
(z, x) \mapsto T(z, x)~; \label{asinc} \feq
any such map is said to be $C^n$ if its derivatives
w.r.t. $z$, $\overline{z}$ and the real variables $x$, up to total order $n$,
exist everywhere and are continuous. The $z$-dependent tensor fields
are acted upon by the derivatives $\partial := \partial/\partial z$,
$\ppartial := \partial/\partial \z$, $\nabla$, $\nabla_{\SS}$ and
their iterates. \parn
\vskip 0.2cm \noindent
\textbf{Iterated derivatives of composite functions.}
Again, we start from the real case. We consider a $u$-dependent
tensor field of arbitrary order $\ell$
\beq H : (-r,r) \times \reali^d \vain \otimes^{\ell} \reali^d~,
\qquad (u,x) \mapsto H(u,x)~. \feq
We associate to $H$ a "tensor" Nemytskij operator sending
a function $f : \reali^d \mapsto (-r,r)$ into the tensor field
\beq H(f,\x) :\reali^d \vain \otimes^{\ell} \reali^d~, \qquad
x \mapsto H(f(x),x)~. \feq
If $H = G = $ a tensor field of order $0$, we recover
definition \rref{defnem} of $G(f,\x)$ and the usual notion of Nemytskij
operator. For $\ell$ arbitrary, we are interested in derivatives of $H(f,\x)$.
Assume $H,f$ to be $C^1$; then
$$ \nabla (H(f,\x))_{\lambda_1 ... \lambda_{\ell+1}} =
\partial_{\lambda_{\ell+1}} \left( H_{\lambda_1 ... \lambda_\ell}(f,\x) \right) = $$
$$ = \left(\partial H_{\lambda_1 ... \lambda_\ell} \right)(f,\x)~
\partial_{\lambda_{\ell+1}} f +
\left(\partial_{\lambda_{\ell+1}} H_{\lambda_1 ... \lambda_\ell}\right)(f,\x) $$
(recall $\partial := \partial/\partial u$). In compact form, this
amounts to write
\beq \nabla \left( H(f,\x) \right) = \partial H (f,\x)
\otimes \nabla f + \nabla H(f,\x)~; \feq
we can symmetrise this result by application of $\SS$ to both
sides, the result being
\beq \nabla_{\SS} \left( H(f,\x) \right) = \partial H (f,\x)
\vee \nabla f + \nabla_{\SS} H(f,\x)~. \label{basic} \feq
Now, we consider the zero order case $H=G : (-r,r) \times \reali^d
\vain \reali$. We are interested in
the iterates $\nabla^m (G(f,\x))$~ $(m \in \naturali_0)$.
These are computed using
as basic tools: the equality $\nabla^m = \nabla_{\SS}^m$ on
real-valued functions, Eq.\rref{basic},
the Leibnitz rule \rref{lei2}, and the fact that $\vee$ is just
the ordinary product when one of the two factors is of order zero. For example,
\beq \nabla \left( G(f,\x) \right) = \nabla_{\SS} (G(f,\x)) =
\partial G(f,\x) \nabla f +
\nabla G (f,\x)~; \label{clear} \feq
\beq \nabla^2 (G(f,\x)) =
\nabla_{\SS} \left(~\partial G(f,\x) \nabla f +
\nabla G(f,\x)~ \right) = \label{clear2} \feq
$$ = \partial^2 G(f,\x) \nabla f \vee \nabla f +
2 \partial \nabla G(f,\x) \vee \nabla f
+ \partial G(f,\x) \nabla^2 f + \nabla^2 G(f,\x)~. $$
To generalize this to any order, we need some
notations. Let us put
\beq \naturali^{(\infty)} := \{ p = (p_1, p_2, ....)~|~p_s \in \naturali
\quad \mbox{for all $s \in \naturali_0$}, \label{ninfty} \feq
$$ \mbox{$p_s \neq 0$ for finitely many $s$}~ \}~; $$
\beq D_{m} := \{ p \in \naturali^{(\infty)}~|~1 \leq p_1 + 2 p_2 + 3 p_3 + ...
\leq m~\} \feq
\beq D_{j w} := \{ p \in \naturali^{(\infty)}~|~p_1 + p_2 + p_3 + ... = j~,
\quad p_1 + 2 p_2 + 3 p_3 + ... = w \} \label{djk} \feq
for all integers $m, j, w$ such that $m \geq 1$, $1 \leq j \leq w$.
The set $D_{j w}$ is non empty for any pair $(j, w)$ as before
({\footnote{In fact, let $\sigma, \rho$ be the
integers such that $\sigma \geq 1$, $ 0 \leq \rho \leq j-1$
and $w = j \sigma + \rho$;
then, $D_{j w}$ contains at least the sequence $p_s :=
(j - \rho) \delta_{s \sigma} +  \rho \delta_{s, \sigma + 1}$.}}), and
$D_m = \cup_{\stackrel{ 1 \leq w \leq m}
{\rule{0mm}{2mm}\scrscr{1 \leq j \leq w}}}~ D_{j w} $
$ = \cup_{\stackrel{ 1 \leq j \leq m}
{\rule{0mm}{2mm}\scrscr{0 \leq \ell \leq m-j}}}~ D_{j, m - \ell}$.
Let us also recall that $\vee^q$ means the $q$-th power w.r.t the
$\vee$ product, see Eq.\rref{veeq}.
\begin{prop}
\label{compos}
\textbf{Proposition.} Let $m \in \naturali_0$; then, there exists a unique
family of positive, integer coefficients
\beq \left( P_{m |p} \right)_{p \in D_m}
\label{family} \feq
such that for each dimension $d$,
each $C^m$ function $G : (-r,r) \times \reali^d
\vain \reali$, $(u,x) \mapsto G(u,x)$ and each $C^m$ function
$f : \reali^d \vain (-r,r)$ it is
\beq \nabla^m \left(G(f,\x)\right) =
\sum_{\stackrel{ 1 \leq j \leq m}
{\rule{0mm}{2mm}\scrscr{0 \leq \ell \leq m-j}}}
\partial^j \nabla^\ell G (f, \x) \vee \label{meloni} \feq
$$ \vee \sum_{p \in D_{j, m - \ell}} P_{m |p}
\left (\vee^{p_1} \nabla^1 f \right)
\vee \left(\vee^{p_2} \nabla^2 f \right) ... ~+
\nabla^m G(f, \x)~. $$
These coefficients are given by
\beq P_{m |p} =
{m! \over (m - p_1 - 2 p_2 - ...)!}~
{1 \over (1!)^{p_1} p_1 ! (2!)^{p_2} (p_2) ! ... }
\label{ee} \feq
for all $p \in D_m$.
\fine
\end{prop}
\textbf{Remarks. i)} For $d=1$ and an $x$-independent $G$, this
statement is known as the Fa\`a di Bruno formula, see e.g. \cite{Com};
we acknowledge S. Paveri Fontana for bibliographical indications
concerning this case. Of course, for $d=1$ the $\vee$ product
collapses into the usual product of real numbers.
Our formulation in terms of tensor products holds for
an arbitrary space dimension $d$ (and $x$-dependent $G$). \parn
Some multidimensional extensions of the Fa\`a di Bruno formula, not using
the tensor language, were proposed in the literature (see, e.g.,
\cite{Gzy} \cite{Mis}); another one was suggested to us by
G. Meloni \cite{Mel}, to whom we are grateful for general discussions
on this subject and its relations with the Stirling numbers. \parn
\textbf{ii}) Eq.\rref{meloni} can be
written in terms of the Bell polynomials and some differential
operators. The Bell polynomials $Y_{w}(t_1, t_2, ... ,t_w)$
\cite{Com}
(which are very similar to the Schur polynomials \cite{Kac})
are defined by the formal
expansion $\mbox{exp}\left(\sum_{r=1}^{\infty} t_r \alpha^r/r! \right) =$
$1 + \sum_{w=1}^{\infty} Y_{w}(t_1, t_2, ..., t_w)~\alpha^w/w!$, and
given explicitly by
\beq Y_{w}(t_1, t_2, ..., t_w) =
\sum_{\stackrel{p \in \naturali^{(\infty)}~,}
{\rule{0mm}{2mm}\scrscr{p_1 + 2 p_2 + 3 p_3 + ... = w}}}
{w! \over (1!)^{p_1} p_1! (2!)^{p_2} p_2! (3!)^{p_3} p_3 ! ...}
{t_1}^{p_1} {t_2}^{p_2} {t_3}^{p_3} ... \feq
for each integer $w \geq 1$. Comparing with Eq.s\rref{meloni}
\rref{ee}, one finds that the formula for the derivatives of a composite
function can be written in the symbolic form
$$ \nabla^m( G(f, \x) ) =
\left(\sum_{0 \leq \ell \leq m-1} \left(\barray{c} m \\ \ell \farray \right)~
Y_{m - \ell}\left(\nabla^1 f \vee \partial, \nabla^2 f
\vee \partial, ...\right) \nabla^\ell \right)(G)(f, \x) + $$
\beq + \nabla^m G(f, \x)~, \feq
the sum denoting an operator which acts on the function
$(u, x) \mapsto G(u, x)$. \fine
\textbf{Proof of Prop \ref{compos}.}
First of all, one can prove by recursion over $m$
the existence of a
family of nonnegative, integer coefficients $(P_{m|p})$
fulfilling \rref{meloni}
for any $d$, $G$ and $f$. The existence of this family
for $m=1$ (and $m=2$) is known from
\rref{clear} (and \rref{clear2}). Assuming existence for some $m$, one
can derive existence at order $m+1$ computing
$\nabla^{m+1} \left(G(f,\x)\right) = \nabla_{S}
\left( \nabla^m \left(G(f,\x)\right) \right)$ from Eq.\rref{meloni}.
Differentiating this equation via the Leibnitz rule \rref{lei2} and
Eq.\rref{basic}, one obtains a sum of terms with nonnegative, integer coefficients,
each one of a type included in the $m+1$ version of
Eq.\rref{meloni}. \parn
Let us pass to prove uniqueness of the coefficients, and derive
the explicit expression \rref{ee}; this will also make clear
that $P_{m|p} \neq 0$ for all $p$.
Let us consider any family $(P_{m|p})$ fulfilling
\rref{meloni} for a given $m$ and all $d$, $G$, $f$. We apply
Eq.\rref{meloni} with this set of
coefficients and $d=1$, calling $\xi$ the $x$ variable; the operator
$\nabla$ is just the ordinary derivative $d/d \xi$, and
\rref{meloni} becomes (adding primes to indices $j$, $\ell$, for future
convenience)
\beq {d^m \over d \xi^m} \left(G(f(\xi),\xi)\right) =
\sum_{\stackrel{ 1 \leq j' \leq m}
{\rule{0mm}{2mm}\scrscr{0 \leq \ell' \leq m-j'}}}
{\partial^{j'+\ell'} G \over \partial u^{j'} \partial \xi^{\ell'}}
(f(\xi), \xi)~ \times \label{melonid} \feq
$$ \times \sum_{p \in D_{j',m-\ell'} } P_{m |p }
\left( {d f(\xi) \over d \xi } \right)^{p_1}
\left( {d^2 f(\xi) \over d \xi^2 } \right)^{p_2} ... +
{\partial^m  G \over \partial \xi^m} (f(\xi), \xi) $$
for all $C^m$ functions $G : (-r, r) \times \reali \vain \reali$
and $f : \reali \vain (-r, r)$. \parn
Now, we
choose a sequence $p = (p_1, p_2, ...)$ in the domain of the
family \rref{family}, and set
\beq j := p_1 + p_2 + ..., \qquad \ell := m - p_1 - 2 p_2 - ... \feq
\beq G(u, \xi) := { u^j \xi^{\ell} \over j! \ell!}~, \qquad
f(\xi) := a_1 \xi + {a_2} {\xi^2 \over 2!} + {a_3} {\xi^3 \over 3!} + ... \feq
where $a_s$ is an arbitrary real parameter if $p_s \neq 0$, and
$a_s := 0$ if $p_s = 0$.
Let us evaluate at $\xi=0$ all derivatives
in \rref{melonid}. The only non vanishing derivatives
$({\partial^{j'+\ell'} G / \partial u^{j'} \partial \xi^{\ell'}})$
at  $(f(0), 0) = (0,0)$ occur for $j' = j$, $\ell' = \ell$ and equal $1$;
furthermore $(d^s f/ d \xi^s)(0) = a_s$. Thus, inserting these choices of $G$, $f$
into Eq.\rref{melonid} and setting $\xi=0$, we obtain
\beq {1 \over j! \ell! }
\left. {d^m \over d \xi^m} \right|_{\xi=0} \left( a_1 \xi
+ a_2 {\xi^2 \over 2!} + ... \right)^j \xi^\ell =
P_{m|p}~ {a_1}^{p_1} {a_2}^{p_2} ...~, \feq
intending $0^0 := 1$; this implies
$$ P_{m|p} = {m! \over j! \ell!} \times \mbox{coefficient of}~
\mbox{${a_1}^{p_1}~{a_2}^{p_2} ... \xi^{m}$
in $\left( \sum_{s} a_s {\xi^s \over s!} \right)^{j}
\xi^{\ell}$}~. $$
Computing the indicated coefficient in the above polynomial,
one obtains for $P_{m|p}$ the expression \rref{ee}; the calculation is
based on the identity
\beq (z_1 + z_2 + z_3 ... )^j =
\sum_{\stackrel{ c_1, c_2, ... \in \naturali~, }
{\rule{0mm}{2mm}\scrscr{ c_1 + c_2 + ... = j  }}}
{j! \over c_1! c_2! ... } {z_1}^{c_1} {z_2}^{c_2} ...~, \feq
to be applied with $z_s := a_s \xi^s / s!$~. \fine
We pass to the complex case, and consider a $z$-dependent
tensor field
\beq H : \B(0,r) \times \reali^d \vain \otimes^{\ell} \complessi^d~,
\qquad (z,x) \mapsto H(z,x)~; \feq
then a tensor Nemytskij operator is defined, sending
a function $f : \reali^d \mapsto \B(0,r)$ into the tensor field
\beq H(f,\x) :\reali^d \vain \otimes^{\ell} \complessi^d~, \qquad
x \mapsto H(f(x),x)~. \feq
Let us discuss the derivatives of $H(f,\x)$: it is readily shown that,
for $H$ and $f$ of class $C^1$,
\beq \nabla \left( H(f,\x) \right) =
\partial H (f,\x) \otimes \nabla f + \ppartial H (f,\x) \otimes
\overline{\nabla f} + \nabla H(f,\x)~, \feq
with $\overline{\nabla f}$ denoting the complex conjugate
of $\nabla f$ (recall Eq.\rref{conj}).
The symmetrised version of this identity is
\beq \nabla_{\SS} \left( H(f,\x) \right) = \partial H (f,\x)
\vee \nabla f + \ppartial H (f,\x) \vee \overline{\nabla f} +
\nabla_{\SS} H(f,\x)~; \label{basicc} \feq
this can be used, together with the Leibnitz rule \rref{lei2} and the
equality $\nabla^m_{\SS} = \nabla^m$ on $\complessi$-valued functions,
to obtain the complex analogue of Prop.\ref{compos}. We formulate
this analogue keeping the notation $\naturali^{(\infty)}$
for the set of sequences \rref{ninfty}, and setting
\beq \D_{m} := \{ (p, q) \in \naturali^{(\infty)}~|~1 \leq p_1 +
2 p_2 + 3 p_3 + ... + q_1 + 2 q_2 + 3 q_3 + ... \leq m~\} \feq
\beq \D_{h k w} :=
\{ (p, q) \in \naturali^{(\infty)} \times
\naturali^{(\infty)}~|~p_1 + p_2 + ... = h~, \quad q_1 + q_2 + ... = k~,
\label{dhk} \feq
$$ p_1 + 2 p_2 + 3 p_3 + ... + q_1 + 2 q_2 + 3 q_3 + ... = w \} $$
for all integers $m, h, k, w$ such that $m \geq 1$, $h, k \geq 0$,
$1 \leq h + k \leq w$. All sets $\D_{h k w}$ are nonempty
({\footnote{The pair $(h, k)$ contains at least
a nonzero element, say $h$.  Let $\sigma, \rho$ be the
integers such that $\sigma \geq 1$, $ 0 \leq \rho \leq h-1$
and $w - k = h \sigma + \rho$;
then, $D_{h k w}$ contains at least the pair $(p, q)$ such that
$p_s := (h - \rho) \delta_{s, \sigma} +  \rho \delta_{s, \sigma + 1}$
and $q_s := k \delta_{s, 1}$.}}), and
$\D_{m} =
\cup_{\stackrel{ 1 \leq w \leq m,~ h, k \geq 0}
{\rule{0mm}{2mm}\scrscr{1 \leq h + k \leq w}}}~~ \D_{h k w} =
\cup_{\stackrel{ h, k \geq 0,~1 \leq h + k \leq m}
{\rule{0mm}{2mm}\scrscr{0 \leq \ell \leq m - h - k}}}~ \D_{h k, m - \ell}$.
\begin{prop}
\label{composc}
\textbf{Proposition.} Let $m \in \naturali_0$.
Then, there exists a unique
family of positive, integer coefficients
\beq \left( P_{m | p q} \right)_{(p, q) \in \D_m} \label{familyc} \feq
such that for each dimension $d$,
each $C^m$ function $G : \B(0,r) \times \reali^d
\vain \complessi$, $(z,x) \mapsto G(z,x)$ and each $C^m$ function
$f : \reali^d \vain \B(0,r)$ it is
\beq \nabla^m \left(G(f,\x)\right) =
\sum_{\stackrel{ h, k \geq 0,~1 \leq h + k \leq m}
{\rule{0mm}{2mm}\scrscr{0 \leq \ell \leq m - h - k}}}
\partial^h \ppartial^k \nabla^\ell G (f, \x)~ \vee \label{melonic} \feq
$$ \vee \sum_{(p,q) \in \D_{h k, m - \ell}} P_{m|p q}
\left (\vee^{p_1} \nabla^{1} f \right)
\vee \left(\vee^{p_2} \nabla^2 f \right) ...
\left (\vee^{q_1} \overline{\nabla^{1} f} \right)
\vee \left(\vee^{q_2} \overline{\nabla^2 f} \right) ...
+ \nabla^m G(f, \x)~. $$
These coefficients are given by
\beq P_{m |p q} =
{m! \over (m - p_1 - 2 p_2 - ... - q_1 - 2 q_2 - ...)!}~\times
\label{eec} \feq
$$ \times~{1 \over (1!)^{p_1} p_1 ! (2!)^{p_2} (p_2) ! ...
(1!)^{q_1} q_1 ! (2!)^{q_2} (q_2) ! ... } $$
for all $(p, q) \in \D_m$.
\end{prop}
\textbf{Proof.} It is similar to the proof
of Prop.\ref{compos}. The existence of a set of nonnegative, integer
coefficients
fulfilling \rref{melonic} for all $d$, $G$ and $f$ is proved by
recursion over $m$. Now, let us consider for given $m$ a family
of coefficients $(P_{m | p q})$
fulfilling \rref{melonic}, and prove these coefficients to be uniquely
determined as in Eq.\rref{eec}, which also implies their positiveness.
To this purpose, it suffices to consider
Eq.\rref{melonic} for $d=1$, calling $\xi$ the $x$ variable
(and adding primes to indices), which gives
\beq {d^m \over d \xi^m} \left(G(f(\xi),\xi)\right) =
\sum_{\stackrel{ h', k' \geq 0, 1 \leq h' + k' \leq m}
{\rule{0mm}{2mm}\scrscr{0 \leq \ell' \leq m- h' - k'}}}
{\partial^{h'+k' + \ell'} G \over
\partial z^{h'} \partial \overline{z}^{k'} \partial \xi^{\ell'}}
(f(\xi), \xi)~
\times \label{melonicd} \feq
$$ \times \sum_{(p, q) \in \D_{h' k', m - \ell'}}
P_{m|p q} \left( {d f(\xi) \over d \xi}
\right)^{p_1} \left( {d^2 f(\xi) \over d \xi^2} \right)^{p_2} ...
\left( \overline{ {d f(\xi) \over d \xi} }\right)^{q_1}
\left( \overline{ {d^2 f(\xi) \over d \xi^2} } \right)^{q_2} ... +
{\partial^m G \over \partial \xi^m}(f(\xi), \xi)~ $$
for all $C^m$ functions $G : \B(0,r) \times \reali \vain
\complessi$ and $f : \reali \vain \B(0, r)$.
\parn
We choose a pair $(p, q)$ in the domain of the family \rref{familyc},
and set
\beq h := p_1 + p_2 + ...~, \quad k := q_1 + q_2 + ....~,
\qquad \ell := m - p_1 - 2 p_2 - ... - q_1 - 2 q_2 - ...~, \feq
\beq G(z, \xi) := { z^h \overline{z}^k \xi^{\ell} \over h! k! \ell!}~, \qquad
f(\xi) := a_1 \xi + {a_2} {\xi^2 \over 2!} + {a_3} {\xi^3 \over 3!} + ... \feq
where $a_s$ is an arbitrary complex parameter if $(p_s, q_s) \neq (0,0)$, and
$a_s := 0$ if $(p_s, q_s) = (0,0)$. Inserting these choices of $G, f$ into
Eq.\rref{melonicd} with $\xi=0$, and working as in the proof of
Prop.\ref{compos}, we readily obtain that $P_{m|p q}$ has the expression
\rref{eec}. \fine
\section{Proof of Lemma \ref{lemma1}. Connections with Prop.s
\ref{compos} and \ref{composc}.}
\label{poly}
We begin with the \parn
\textbf{Proof of Lemma \ref{lemma1} }. Let us
first prove the uniqueness of the polynomial fulfilling
\rref{dmxm} for all functions $F$.
To this purpose we consider any polynomial $P$ fulfilling
\rref{dmxm}, and show that $P((\nu_{h \ell}), \rho)$
is uniquely determined for all $\nu_{h \ell} \in \reali$
and $\rho \in (0,+\infty)$, a fact yielding the thesis. To prove
this fact, given $(\nu_{h \ell})$, $\rho$ as above we
define a function
\beq F = F_{m, (\nu_{h \ell}), \rho} : \reali^2 \vain
\reali~, \feq
$$ F_{m, (\nu_{h \ell}), \rho}(u,\xi) :=
 \sum_{\stackrel{ 1 \leq h \leq m}
{\rule{0mm}{2mm}\scrscr{0 \leq \ell \leq m - h}}}
\nu_{h \ell}~
{(u - \rho)^h (\xi - \ln \rho)^\ell \over h! \ell!}~. $$
Then, it is found that
\beq {\partial^{h+\ell} F_{m, (\nu_{h \ell}), \rho}
\over \partial u^h \partial \xi^\ell}(\rho,
\ln \rho ) = \nu_{h \ell}~, \qquad
{\partial^{m} F_{m, (\nu_{h \ell}), \rho}
\over \partial \xi^m}(\rho, \ln \rho) = 0~; \feq
these results, with Eq.\rref{dmxm} for $P$, imply
\beq P((\nu_{h \ell}), \rho) =
\left. {d^m \over d \xi^m} \right|_{\xi=\ln \rho}~
F_{m, (\nu_{h \ell}), \rho}(e^\xi, \xi)~, \label{imp} \feq
and so $P$ is uniquely determined. \parn
The existence of the polynomial $P = P_m$ fulfilling
Eq.\rref{dmxm}, its form
\rref{form} and Eq.\rref{rec}
can be proved in a single step.
In fact, it is easy to check that the polynomial
$$ P_m((\nu_{j \ell}), \rho) :=
\sum_{\stackrel{ 1 \leq j \leq m}
{\rule{0mm}{2mm}\scrscr{0 \leq \ell \leq m-j}}} P_{m j \ell}~
\nu_{j \ell}~ \rho^j~ $$
fulfills \rref{dmxm} for each $m \in \naturali_0$, if
the coefficients $(P_{m j \ell})$ satisfy the recursion
relations \rref{rec}.
Let us consider these coefficients more closely. From
the recursion relations, it is clear that they are integers.
By construction, we have
\beq P_{m j \ell} = P_m \Big( \nu_{j'  \ell'} = \delta_{j' j}~
\delta_{\ell' \ell}~, \rho = 1 \Big)~; \feq
the r.h.s. is given by Eq.\rref{imp} for $P = P_m$, considering
the derivative at $\xi=0$ and the
function $F_{m, (\nu_{j' \ell'} = \delta_{j' j}
\delta_{\ell' \ell}), \rho=1}(u, \xi) = (u - 1)^j \xi^\ell/(j! \ell!)$.
In conclusion
$$ P_{m j \ell}  =
\left. {d^m \over d \xi^m} \right|_{\xi=0}~{(e^\xi - 1)^j \xi^\ell
\over j! \ell!}~, $$
and this proves Eq.\rref{equi1}. Finally, let us write
$(e^\xi - 1)^j = \sum_{r=1}^j \left( \barray{c} j
\\ r \farray \right) e^{r \xi} (-1)^{j - r}$ and expand $e^{r \xi}$
about $\xi=0$; inserting the expansion in the previous equation, and isolating
the coefficient of $\xi^m$ we obtain the expression \rref{equi2} for
$P_{m j \ell}$. \fine
Now, let us point out a relation between the coefficients
in the general formulae for the derivation of a composite function and
the coefficients $P_{m j \ell}$ of the polynomials $P_m$ (an equivalent
formulation of the forthcoming Eq.\rref{rel1}
for $\ell=0$ can be found in \cite{Com}).
\begin{prop}
\label{coeff}
\textbf{Lemma.} Consider the coefficients $P_{m |p}$
and $P_{m h k \ell|p q}$ in Eq.s\rref{meloni}
\rref{melonic} for the derivation of a composite function, in the
real and complex cases respectively. These are related to the
$P_{m j \ell}$ coefficients of the polynomials $P_m$ by
\beq \sum_{p \in D_{j, m - \ell} }
P_{m | p} = P_{m j \ell}~; \label{rel1} \feq
\beq \sum_{(p, q) \in \D_{h k, m - \ell}} P_{m |p q}~ =
\left( \barray{c} h + k \\ h
\farray \right) P_{m, h + k, \ell}~. \label{rel2} \feq
\end{prop}
\textbf{Proof.} We will derive in more detail Eq.\rref{rel2}; the argument
yielding to \rref{rel1} is even simpler. Consider any $C^m$ function
\beq G : \complessi \times \reali \vain \reali~,
\qquad (z, \xi) \vain G(z, \xi)~, \feq
and apply the "$1$-dimensional"
formula \rref{melonicd} with
$f : \reali \vain \reali$, $f(\xi) := e^\xi$; this gives
\beq {d^m \over d \xi^m} \left(G(e^\xi,\xi)\right) = \label{theprev} \feq
$$ = \sum_{\stackrel{ h', k' \geq 0, 1 \leq h' + k' \leq m}
{\rule{0mm}{2mm}\scrscr{0 \leq \ell \leq m - h' - k'}}}
{\partial^{h' + k' + \ell'} G \over \partial z^{h'}
\partial \overline{z}^{k'} \partial \xi^{\ell'}} (e^\xi, \xi)~
\Big( \sum_{(p, q) \in \D_{h' k', m - \ell'} } P_{m | p q} \Big)~
e^{(h' + k') \xi} + {\partial^m  G \over \partial \xi^m} (e^\xi, \xi) ~. $$
Let us choose, in particular,
\beq G(z, \xi) := { (z - 1)^{h} (\overline{z} - 1)^{k} \xi^{\ell} \over
h! k! \ell! } \feq
for fixed $h, k, \ell \geq 0$ with $1 \leq h + k \leq m$ and
$0 \leq \ell \leq m - h - k$; we apply Eq.\rref{theprev}
with this $G$ at $\xi=0$. The l.h.s. of \rref{theprev} gives
the derivative at $\xi=0$ of $G(e^\xi, \xi) =
(e^\xi - 1)^{h + k} \xi^{\ell} / (h! k! \ell!) $; in the r.h.s.
the only nonzero derivative is
$(\partial^{h + k + \ell} G / \partial z^h
\partial \overline{z}^{k} \partial \xi^{\ell}) (1, 0) = 1$, and thus
\beq \left. {d^m \over d \xi^m} \right|_{\xi=0}
{ (e^\xi - 1)^{h + k} \xi^{\ell} \over {h}! {k}! {\ell!} } =
\sum_{(p, q) \in \D_{h k, m - \ell} } P_{m | p q}~. \feq
By comparison with Eq.\rref{equi1}, we see that the l.h.s. in this
equation equals $\left( \barray{c} h + k \\ h
\farray \right)$ $ P_{m, h + k, \ell}$; so, Eq.\rref{rel2} is proved. \parn
Eq.\rref{rel1} is derived similarly, applying formula
\rref{melonic} to $f(\xi) := e^\xi$ and $G : \reali^2 \vain \reali$,
$G(u, \xi) := (u - 1)^{j} \xi^\ell / (j! \ell!)$.
\fine
\section{More on tensor functional spaces. Some inequalities.}
\label{ineq}
We always consider spaces of complex tensor-valued functions;
the real analogues can be regarded as subspaces of the complex
ones. \vskip 0.2cm \noindent
\textbf{H\"older type inequalities.} $L^p$ spaces of functions
$ T : \reali^d \mapsto \otimes^\ell \complessi^d$ have been already
introduced, with the norms $\| T \|_{L^p} := \|~ | T |~ \|_{L^p}$
(see Eq.s (\ref{nt}-\ref{nlpf})~). \parn
Let $h \in \naturali_0$, $\ell_1, ..., \ell_h \in \naturali$
and $p_1, ..., p_h, r \in [1,+\infty]$ be such that
$1/p_1 + ... + 1/p_h = 1/r$; further, consider the tensor-valued functions
$T_1 \in L^{p_1}(\reali^d, \otimes^{\ell_1} \complessi)$, ...,
$T_h \in L^{p_h}(\reali^d, \otimes^{\ell_h} \complessi)$. Then
$T_1 \otimes .... \otimes T_h$, $T_1 \vee .... \vee T_h
 \in L^r(\reali^d, \otimes^{\ell_1 + ... + \ell_h}
\complessi^d)$, and
\beq \| T_1 \otimes ... \otimes T_h \|_{L^r}~,~
\| T_1 \vee ... \vee T_h \|_{L^r} \leq \| T_1 \|_{L^{p_1}} ...
\| T_h \|_{L^{p_h}}~; \feq
these H\"older type inequalities follow readily
from the classical H\"older inequality and from
Eq.\rref{notv} for the tensor norms $|~|$. \parn
\vskip 0.2cm \noindent
\textbf{Fourier transform and Sobolev spaces.}
We denote with $\FF : \SS'(\reali^d, \complessi)
\vain \SS'(\reali^d, \complessi)$ the Fourier transform
for tempered distributions;
when using this, we write $k = (k_1, ..., k_d)$
for the running variable on $\reali^d$ and $\k : \reali^d \vain \reali^d$
for the identity mapping of this space into itself. We normalise
$\FF$ so that $(\FF f)(k) = (2 \pi)^{-d/2} \int_{\reali^d}
d x ~e^{- i k {\scriptstyle{\bullet} } x} f(x)$ for $f$ in $L^1$. \parn
For $m \in \naturali$, we denote with $\SS'(\reali^d, \otimes^m \complessi^d)$
the space of families $T = (T_{\lambda})_{\lambda \in \{1, ..., d\}^m}$
where each component is in
$\SS'(\reali^d, \complessi)$. We extend
componentwisely the Fourier transform to tensor valued
tempered distributions, setting
\beq \FF T := (\FF T_{\lambda})_{\lambda \in \{1,...,d\}^m} \feq
for each $T$ as before.
For $f \in \SS'(\reali^d, \complessi)$ and $m \in \naturali$,
the Fourier representation of distributional derivatives can be written as
\beq \nabla^m f = \FF^{-1}( \otimes^{m} (i \k) \FF f )~, \label{nabm} \feq
where $\otimes^m(i \k) := (i \k) \otimes ... \otimes (i \k)$ ($m$
times). This fact, with the preservation of the $L^2$ norm under
$\FF^{-1}$, implies
\beq \| \nabla^m f \|_{L^2} = \| \otimes^m (i \k) \FF f \|_{L^2} =
\| | \k |^m \FF f \|_{L^2}~, \label{pres} \feq
whenever one of the above three functions is in $L^2$; the second
equality follows from the fact that the tensor norm
$| \otimes^m  (i k) |$ equals $| k |^m = \sqrt{k_1^2 + ... + k_d^2}^{~m}$
for each $k \in \reali^d$ (recall Eq.\rref{notv}).
Due to these facts and to the equality
\beq (1 + | \k |^2)^n = \sum_{m=0}^n
\left( \barray{c} n \\ m \farray \right) | \k |^{2 m}~, \feq
we see from Eq.s\rref{hnab} \rref{nonab} that
\beq H^n(\reali^d, \complessi)
= \{ f \in \SS'(\reali^d, \complessi)~|~
\sqrt{1 + | \k |^2}^{~n} \FF f \in L^2(\reali^d, \complessi)~\}~,
\label{hnabf} \feq
\beq \| f \|_n = \| \sqrt{1 + | \k |^2}^{~n} \FF f \|_{L^2}~; \label{nonabf} \feq
thus, the conventional definition of Bessel potential spaces
is recovered (Eq.s\rref{hnabf}
\rref{nonabf} also make sense for non integer $n$, and
could be written in terms of the operator $\sqrt{1 - \Delta}^{~n} :=
\FF^{-1} \sqrt{1 + | \k |^2}^{~n} \FF$).
\vskip 0.2cm\noindent
\textbf{Hausdorff-Young inequality.} We recall the definition
\rref{defe} for the function $E(~)$.
Let $p \in [1,2]$, $r \in [2,+\infty]$ be such that $1/p+1/r = 1$.
Further, let $m \in \naturali$; for all $T \in L^p(\reali^d,
\otimes^m \complessi^d)$, it is $\FF^{-1} T \in L^{r}(\reali^d,
\otimes^m \complessi)$
and
\beq \| \FF^{-1} T \|_{L^r} \leq C_{r d}~ \| T \|_{L^p}~, \qquad
C_{r d} := {1 \over (2 \pi)^{d/2 - d/r}}~
\left( {E(1/r) \over E(1 - 1/r)} \right)^{d/2}~; \label{crd} \feq
furthermore, the above constant is the sharp one for the written
inequality. In the case $m=0$, where $T$ is $\complessi$-valued,
this result is classic (see Beckner's
paper \cite{Bec}, or \cite{Lie2}; these give a slightly
different expression for $C_{r d}$ due to a different normalisation
for the Fourier transform). For $m=1,2,3,...$, the above statement
follows from the Marcinkiewicz-Zygmund theorem on the componentwise
vector extensions of operators between $L^p$ spaces (see, e.g.,
\cite{Gar}, Ch.5, Theor.2.7; it must be noted that the
norm \rref{nt} we use on $\otimes^{m} \complessi^d$ is a Hilbertian,
$\ell^2$-type norm, as required by the theorem. We acknowledge
L. Colzani for pointing out to us the cited result, and W. Beckner for
general indications on the vectorial extensions of the Hausdorff-Young
inequality).
\vskip 0.2cm \noindent
\textbf{An interpolation inequality.} We rephrase in tensor
notations a well known inequality; the aim is to point out
the absence of $d$-dependent constants.
\begin{prop}
\textbf{Proposition.} Let $m \in \naturali$,
$f \in L^2(\reali^d, \complessi)$ such that
$\nabla^m f \in L^2(\reali^d, \otimes^m \complessi)$.
Then, for $\ell \in \{0,...,m\}$,
it is $\nabla^\ell f \in L^2(\reali^d, \otimes^\ell \complessi)$ and
\beq \| \nabla^\ell f \|_{L^2} \leq \| f \|^{1 - \ell/m}_{L^2}~
\| \nabla^m f \|^{\ell/m}_{L^2}~\leq (1 - {\ell\over m}) \| f \|_{L^2} +
{\ell \over m} \| \nabla^m f \|_{L^2} \label{inter} \feq
(intending $\ell/m := 1$ for $\ell = m = 0$ and $0^t := 0$ for each $t$).
\end{prop}
\textbf{Proof.}
Let us prove that $\nabla^\ell f$ is $L^2$
and satisfies the first inequality \rref{inter}.
Due to \rref{pres}, it suffices to show that
$| \k |^\ell \FF f \in L^2(\reali^d, \complessi)$, and
\beq \|~ | \k |^\ell \FF f  \|_{L^2} \leq \| \FF f \|^{1 - \ell/m}_{L^2}~
\|~ | \k |^m \FF f \|^{\ell/m}_{L^2}~; \feq
the desired result follows writing
$| \k |^\ell | \FF f | =$ $\left( | \FF f |^{1 - \ell/m}
\right)  \left( | \k |^\ell | \FF f |^{\ell/m} \right) $
and using the H\"older inequality $\|~\|_{L^2} \leq \|~ \|_{L^p}
\| ~ \|_{L^q}$, with $1/p = 1/2 - \ell/(2 m)$,
$1/q = \ell/(2 m)$. \parn
The second inequality \rref{inter} follows trivially from the
first one, due to the Young inequality $v^{1 - \theta} w^{\theta}
\leq (1 - \theta) v + \theta w$ ($v, w \in [0, +\infty)$,
$\theta \in [0,1]$).
\fine
\vskip 0.2cm \noindent
\textbf{A Gagliardo type inequality.} The structure of the inequality
is familiar, but little seems available about the
constant appearing therein. Here we propose an explicit (probably
non sharp) estimate for the constant, based on the use
of Fourier methods and the Hausdorff-Young inequality \rref{crd}.
All derivatives in the forthcoming statement are intended in the
distributional sense. \parn
\begin{prop}
\label{gagl}
\textbf{Proposition.} Let $\ell, m, a \in \naturali$~, $\ell \leq m$~,
$a > {d / 2}$; intend $m/\ell :=1$ if $\ell=m=0$.
If $f \in H^a(\reali^d, \complessi)$
and $\nabla^m f \in L^2(\reali^d, \otimes^m \complessi^d)$, then
$\nabla^\ell f \in L^{2 m/\ell}(\reali^ d, \otimes^\ell \complessi)$ and
\beq \| \nabla^\ell f \|_{L^{2 m/\ell}} \leq
\left( E(\ell/(2 m)) \over E(1 - \ell/(2 m)) \right)^{d/2}
\left(S_{a d}~ \| f \|_a \right)^{1 - \ell/m}\| \nabla^m f \|^{\ell/m}_{L^2}~,
\label{tesi} \feq
where $E(~)$ is defined as in \rref{defe} and $S_{a d}$ is the constant
\rref{sad}.
\end{prop}
\textbf{Proof.} We put
\beq r := {2 m \over \ell}~, \qquad p := {2 m \over 2 m - \ell}~,
\qquad s := {2 m \over m - \ell} \feq
(intending $s := \infty$ if $m=\ell \neq 0$;
$r := s:= 2$ and $p:= 1$ if $m=\ell=0$).
Then $r,s \in [2, +\infty]$, $p \in [1,2]$ and
\beq {1 \over p} + {1 \over r} = 1~, \qquad
{1 \over r} + {1 \over s} = {1 \over 2}~,
\qquad {2 \over s} + {1 \over r} = {1 \over p}~. \feq
By Eq.\rref{nabm} and the Hausdorff-Young
inequality \rref{crd}, we have
\beq \| \nabla^{\ell} f \|_{L^r} =
\| \FF^{-1}\left( \otimes^{\ell}(i \k) \FF f \right) \|_{L^r} \leq
C_{r d}~ \| \otimes^{\ell}(i \k) \FF f \|_{L^p} = \label{muno} \feq
$$ = {1 \over (2 \pi)^{d/s} }~
\left( {E(1/r) \over E(1 - 1/r)} \right)^{d/2}
\| | \k |^{\ell} \FF f \|_{L^p}~,
$$
provided that the last function be $L^p$. To prove this and the
rest, we write
\beq | \k |^{\ell} | \FF f | = {1 \over (1 + | \k |^2)^{a/s}}~
\left( (1 + | \k |^2 )^{a/s} | \FF f |^{2/s} \right)~
\left( | \k |^{\ell} | \FF f |^{2/r} \right) \feq
and apply the H\"older inequality to the above three factors, recalling
that $1/s + 1/s + 1/r=1/p$; this gives
\beq \|~ | \k |^{\ell} | \FF f | \|_{L^p}
\leq \| {1 \over (1 + | \k |^2)^{a/s}} \|_{L^s}~
\| (1 + | \k |^2)^{a/s}  | \FF f |^{2/s} \|_{L^s}~
\| | \k |^{\ell} | \FF f |^{2/r} \|_{L^r}~. \label{zero} \feq
Let us show that the norms of the three factors are finite, and
compute them. We have
\beq \| {1 \over (1 + | \k |^2)^{a/s}} \|_{L^s} =
\left( \int_{\reali^d} {d k \over (1 + | k |^2)^a} \right)^{1/s}
= (2 \pi)^{d/s} S_{a d}^{2/s}~; \label{uno} \feq
\beq \| (1 + | \k |^2)^{a/s}  | \FF f |^{2/s} \|_{L^s} =
\left( \int_{\reali^d} d k (1 + | k |^2)^a | \FF f |^{2} \right)^{1/s}
= \| f \|_{a}^{2/s}~; \label{due} \feq
\beq \| | \k |^{\ell} | \FF f |^{2/r} \|_{L^r} =
\left( \int_{\reali^d} d k | k |^{\ell r} | \FF f |^{2} \right)^{1/r} =
\left( \int_{\reali^d} d k | k |^{2 m} | \FF f |^{2} \right)^{1/r} =
\label{tre} \feq
$$  =\| \otimes^{m} (i \k) \FF f \|_{L^2}^{2/r} =
\| \nabla^m f \|_{L^2}^{2/r}~ $$
(some intermediate steps where integrals appear do not make sense for
$s$ or $r = \infty$, but the final results also cover these cases). \parn
Now, we insert Eq.s (\ref{uno}-\ref{tre}) into \rref{zero}, and the
result into Eq.\rref{muno}; after explicitating $r$ and $s$, we readily
get the thesis \rref{tesi}. \fine
\vskip 0.2cm \noindent
\textbf{The Adams-Frazier inequality.} Again, we discuss an inequality whose
structure is known \cite{Ada2} but for which the constants
were not previously estimated, to the best of our knowledge. As
in the previous subsection, we intend derivatives in the distributional sense.
\begin{prop}
\label{adam}
\textbf{Proposition.} Let $h, k \in \naturali$, $h+k \geq 1$
and $i_1, ..., i_h, g_1, ..., g_k$ $\in \naturali_0$; put
$m := i_1 + ... + i_h + g_1 + ... + g_k$. Furthermore, let
$a > d/2$ and consider a function
$f \in H^a(\reali^d, \complessi)$ with
$\nabla^m f \in L^2(\reali^d, \complessi)$. Then
$\nabla^{i_1} f , ..., \nabla^{g_k} f$ are ordinary functions,
the $\vee$ product written below is $L^2$, and
\beq \| \nabla^{i_1} f \vee ... \vee \nabla^{i_h} f
\vee \overline{\nabla^{g_1} f} ... \vee \overline{\nabla^{g_k} f} \|_{L^2}
\leq \label{adin} \feq
$$ \leq U_{m, h + k, d} (S_{a d } \| f \|_a)^{h+k-1} \| \nabla^m f \|_{L^2}~,
$$
where $U_{m, h + k, d}$~$(\leq 1)$ is defined as in \rref{defu}.
\end{prop}
\textbf{Proof.} The statement is trivial if $h+k=1$; in the sequel
$h+k \geq 2$. We put
\beq p_u := {2 m \over i_u}~, \qquad q_v := {2 m \over g_v} \feq
for $u \in \{1, ..., h\}$, $v \in \{1,..., k\}$; then $p_u, q_v \in [2,+\infty)$
and
\beq \sum_{u=1}^h {1 \over p_u} + \sum_{v=1}^k {1 \over q_v} = {1 \over 2}~.
\feq
Thus, from the H\"older inequality we get
\beq \| \nabla^{i_1} f \vee ...
\vee \overline{\nabla^{g_k} f} \|_{L^2} \leq
\| \nabla^{i_1} f \|_{L^{p_1}}~
... \| \overline{\nabla^{g_k} f} \|_{L^{q_k}}
\label{hold} \feq
provided that all norms in the r.h.s. be finite.
In fact it is so, and the Gagliardo inequality
with $\ell=i_u$ or $\ell=g_v$ and $m$ as above gives
\beq \| \nabla^{i_u} f \|_{L^{p_u}} \leq
\left( {E(i_u/(2 m)) \over E(1 - i_u/(2 m))} \right)^{d/2}
\left(S_{a d} \| f \| \right)^{1 - i_u/m}_a
\| \nabla^m f \|^{i_u/m}_{L^2}~, \feq
$$ \| \overline{\nabla^{g_v} f} \|_{L^{q_v}}
= \| \nabla^{g_v} f \|_{L^{q_v}} \leq
\left( {E(g_v/(2 m)) \over E(1 - g_v/(2 m))} \right)^{d/2}
\left( S_{a d} \| f \| \right)^{1 - g_v/m}_a \| \nabla^m f \|^{g_v/m}_{L^2}~. $$
Inserting these estimates into \rref{hold} and recalling that
$i_1 + ... + g_k = m$, we readily get
\beq \| \nabla^{i_1} f \vee ...
\vee \overline{\nabla^{g_k} f} \|_{L^2} \leq
\UU_{i_1 ... g_k, d} (S_{a d } \| f \|_a)^{h+k-1} \| \nabla^m f \|_{L^2}
~, \feq
\beq \UU_{i_1 ... g_k, d} := \left(
\Pi_{u=1}^h {E(i_u/(2 m)) \over E(1 - i_u/(2 m))}~
\Pi_{v=1}^k {E(g_v/(2 m)) \over E(1 - g_v/(2 m))}~
\right)^{d/2}~; \feq
now, the proof is concluded if we show that
\beq \UU_{i_1 ... g_k, d} \leq U_{m, h+k, d}~. \label{ifwe} \feq
To prove this, we put
\beq \HH : (0,1)^{h+k} \vain \reali~, \feq
$$ \vartheta = (\vartheta_1, ..., \vartheta_{h+k}) \mapsto
\HH(\vartheta) := \sum_{w=1}^{h+k}
{\vartheta_w \over 2} \ln {\vartheta_w \over 2} -
(1 - {\vartheta_w \over 2}) \ln (1 - {\vartheta_w \over 2})
$$
and denote with $\eta$, $\zeta_{(1)}$ ,..., $\zeta_{(h+k)}$
the following points of the cube $(0,1)^{h+k}$:
\beq \eta := ({i_1 \over m}, ... , {i_h \over m},
{g_1 \over m}, ... ,{g_k \over m})~,
\qquad \zeta_{(1)} := (1 - {h+k-1 \over m}, {1 \over m}, ..., {1 \over m})~, \feq
$$ \zeta_{(2)} :=
({1 \over m}, 1 - {h+k-1 \over m}, {1 \over m}, ..., {1 \over m})~,
~~...~~,~\zeta_{(h+k)} :=
({1 \over m}, .... , {1 \over m}, 1 - {h+k-1 \over m})~; $$
then
\beq {d \over 2}~ \HH(\eta) = \ln \UU_{i_1 ... g_k, d}~, \qquad
{d \over 2}~ \HH(\zeta_{(1)}) = ... = {d \over 2}~ \HH(\zeta_{(h+k)}) =
\ln U_{m, h + k, d}~. \label{zeta} \feq
It is readily cheched that $\HH$ is a convex function (its Hessian is a
diagonal, positive defined matrix); on the other hand
\beq \eta = \sum_{w=1}^{h+k} t_{w} \zeta_{(w)}~, \feq
$$ t_1 := {i_1 - 1\over m - h - k},~ ...,~
t_{h} := {i_h - 1\over m - h - k},~
t_{h+1} := {g_1 - 1 \over m - h - k},~ ...,~
t_{h+k} := {g_k - 1 \over m - h - k} $$
(intending $t_1 := ... = t_{h+k}:= 1/(h+k)$ if $m=h+k$;
this case occurs if and only if $i_1 = ... = g_k = 1$). By construction
$t_w \geq 0$ for all $w$, and $\sum_{w=1}^{h+k} t_{w} = 1$. Thus, by
the convexity of $\HH$ and \rref{zeta},
\beq \ln \UU_{i_1 ... g_k, d} = {d \over 2}~ \HH(\eta) \leq
{d \over 2}~\sum_{w=1}^{h+k} t_w \HH(\zeta_{(w)}) = \ln U_{m, h + k, d}~. \feq
This gives Eq.\rref{ifwe}, and concludes the proof. \fine
\section{Proofs of Propositions \ref{nemreale} and \ref{nemc}}
\label{last}
We work out in detail the complex case only, corresponding to Prop.\ref{nemc};
the real case
can be treated similarly (and is in fact simpler). \parn
The main step towards the proof of Prop.\ref{nemc} is to estimate
the norms $\| \nabla^m (G(f, \x)) - \nabla^m(G(0, \x)) \|_{L^2}$
$(0 \leq m \leq n)$ under the given assumptions on $G$ and $f$.
Formula \rref{melonic} for the derivatives of a composite function is
essential for this purpose; the result is the following.
\vskip 0.2cm \noindent
\begin{prop}
\label{dmg}
\textbf{Proposition.} Let $n, a \in \naturali$ and $a > d/2$.
Consider a function
$G : \B(0,r) \times \reali^d \vain \complessi$ with the $\Pi^{n}$ property,
such that $G(0,\x) \in H^{n}(\reali^d, \complessi)$,
and a function $f \in H^{n}(\reali^d, \complessi) \cap
H^a_{r}(\reali^d, \reali)$. \parn
Then $G(f,\x) \in H^{n}(\reali^d, \complessi)$; furthermore, it is
\beq \| G(f, \x) - G(0, \x) \|_{L^2} \leq
\flat_{0}(G, S_{a d} \| f \|_a ) \| f \|_{L^2}~,
\label{zget} \feq
\beq \| \nabla^m \left(G(f,\x)\right)
- \nabla^m \left(G(0,\x)\right) \|_{L^2} \leq
\beta_{m d}(G, S_{a d} \| f \|_a)~ \| \nabla^m f \|_{L^2} + \label{wget} \feq
$$ + \left( b_{m d}(G, S_{a d} \| f \|_a)
+ \flat_m(G, S_{a d} \| f \|_a) \right)~ \| f \|_{L^2} \qquad\qquad (1 \leq
m \leq n)~, $$
where $\flat_0$, $\flat_m$, $\beta_{m d}$, $b_{m d}$ are defined as
in Eq.s\rref{sim}, \rref{betam} and \rref{bm}.
\end{prop}
\textbf{Proof.}
\textsl{Step 1: $G(f,\x) \in L^2(\reali^d, \complessi)$, and
Eq.\rref{zget} holds}. By our assumptions $G(0, \x)$ is $L^2$; thus,
we must show that $G(f, \x) - G(0, \x)$ is $L^2$,
and that its norm has the bound \rref{zget}.
Recalling the definition \rref{simc} we see that, for $\rho \in [0, r)$,
$z \in \Bc(0, \rho)$ and $x \in \reali^d$,
\beq | G(z, x) - G(0, x) |
\leq \flat_{0}(G, \rho) | z |~. \label{ze} \feq
Due to \rref{imbed}, we can apply this with $z = f(x)$,
$\rho = S_{a d} \| f \|_a$;
thus, integrating over $x$ (the squares of) both sides, we obtain
Eq.\rref{zget}.
\vskip 0.1cm \noindent
\textsl{Step 2 (for $n \geq 1$).
It is $\nabla^m \left(G(f,\x)\right) \in L^2(\reali^d, \complessi)$
and \rref{wget} holds for $1 \leq m \leq n$, if $f \in
H^{n}(\reali^d, \complessi) \cap
H^a_{r}(\reali^d, \reali) \cap C^n(\reali^d, \complessi)$.} The
supplementary assumption that $f$ be $C^n$ allows us to use Eq.\rref{melonic} for
the derivatives of a composite function, in the usual (i.e., not distributional)
sense (next, the $C^n$ assumption on $f$ will be removed
by appropriate density arguments, see Step 3). For the sake of brevity,
let us put
\beq \GG : \B(0, r) \times \reali^d \vain \complessi~, \qquad
\GG(z, x) := G(z, x) - G(0,x)~; \feq
then $\GG$ also possesses the $\Pi^{n}$ property.
By our assumptions, $\nabla^m(G(0, \x))$ is $L^2$; we will
show that $\nabla^m(\GG(f, \x))$ is also $L^2$, with its $L^2$ norm
bounded by the r.h.s. of \rref{wget}.
Let us express $\nabla^m (\GG(f, \x))$ according to
Eq.\rref{melonic} (with $G$ replaced by $\GG$); this implies
\beq \| \nabla^m \left(\GG(f,\x)\right) \|_{L^2} \leq
\sum_{\stackrel{ h, k \geq 0,~1 \leq h + k \leq m}
{\rule{0mm}{2mm}\scrscr{0 \leq \ell \leq m - h - k}}}
\| \partial^h \ppartial^k \nabla^\ell \GG (f, \x) \|_{L^{\infty}} \times
\label{melnorm} \feq
$$ \times \sum_{(p,q) \in \D_{h k, m - \ell}} P_{m|p q}
\| \left(\vee^{p_1} \nabla^{1} f \right)
\vee \left(\vee^{p_2} \nabla^2 f \right) ...
\left (\vee^{q_1} \overline{\nabla^{1} f} \right)
\vee \left(\vee^{q_2} \overline{\nabla^2 f} \right) ... \|_{L^2}
+ \| \nabla^m \GG(f, \x) \|_{L^2}~, $$
provided that all norms in the r.h.s. be finite (we have used the
H\"older inequality for tensor fields, in the form
$\| T \vee S \|_{L^2} \leq \| T \|_{L^{\infty}} \| S \|_{L^2}$). \parn
Let $h, k, \ell$ be as in the above sum;
then $\partial^h \ppartial^k \nabla^\ell \GG =
\partial^h \ppartial^k \nabla^\ell G$ (because the
$z$-independent term $x \mapsto G(0, x)$
disappears on application of $\partial^h \ppartial^k$). Recalling
the definition \rref{sihkc} we see that, for $\rho \in [0, r)$,
$z \in \Bc(0, \rho)$ and $x \in \reali^d$,
$$ | \partial^h \ppartial^k \nabla^\ell \GG(z, x) | =
| \partial^h \ppartial^k \nabla^\ell G(z, x) | \leq
\sigma_{h k \ell}(G, \rho)~; $$
due to \rref{imbed}, we can apply this inequality with $z = f(x)$ and $\rho =
S_{a d} \| f \|_a$, so
\beq \| \partial^h \ppartial^k \nabla^\ell \GG(f, \x) \|_{L^{\infty}} =
\sup_{x \in \reali^d}
| \partial^h \ppartial^k \nabla^\ell G(f(x), x) | \leq
\sigma_{h k \ell}(G, S_{a d} \| f \|_a)~. \label{rsup} \feq
Let us pass to the products
$$ \left(\vee^{p_1} \nabla^{1} f \right)
\vee \left(\vee^{p_2} \nabla^2 f \right) ...
\left (\vee^{q_1} \overline{\nabla^{1} f} \right)
\vee \left(\vee^{q_2} \overline{\nabla^2 f} \right) ... ~$$
with $(p, q) \in \D_{h k,m - \ell}$; any one of them can be reexpressed as
$$ \nabla^{i_1} f \vee ... \nabla^{i_j} f \vee
\overline{\nabla^{g_1} f} \vee ... \vee \overline{\nabla^{g_k} f}  $$
where $i_1 \leq i_2 \leq ... \leq i_j, g_1 \leq g_2 ... \leq g_k$
are positive integers, determined in an obvious way by the sequences
$(p, q)$. The above expression envolves $h+k$ factors, with
total order $i_1 + ... + g_k = m - \ell$, and its $L^2$ norm
can be bounded by the Adams-Frazier inequality \rref{adin} in terms of
$\| \nabla^{m - \ell} f \|_{L^2}$; the last norm can be interpolated
with the second inequality \rref{inter}. Summing up
$$ \| \left(\vee^{p_1} \nabla^{1} f \right)
\vee \left(\vee^{p_2} \nabla^2 f \right) ...
\left (\vee^{p_1} \overline{\nabla^{1} f} \right)
\vee \left(\vee^{p_2} \overline{\nabla^2 f} \right) ... \|_{L^2} = $$
$$ = \| \nabla^{i_1} f \vee ... \vee \overline{\nabla^{g_k} f} \|_{L^2}
\leq U_{m - \ell, h + k, d}
(S_{a d } \| f \|_a)^{h+k-1} \| \nabla^{m - \ell} f \|_{L^2} \leq  $$
\beq
\leq U_{m - \ell, h + k, d} (S_{a d } \| f \|_a)^{h+k-1}
\left( {\ell \over m} \| f \|_{L^2} +
(1 - {\ell \over m}) \| \nabla^m f \|_{L^2} \right)~.
\label{rvee} \feq
Finally, let us consider the term $\| \nabla^m \GG(f, \x) \|_{L^2}$;
a straightforward generalization of Eq.\rref{ze} implies
\beq \| \nabla^{m} \GG(f, \x) \|_{L^2} \leq
\flat_{m}(G, S_{a d} \| f \|_a ) \| f \|_{L^2}~.
\label{rget} \feq
Let us insert the results \rref{rsup} \rref{rvee} \rref{rget} into
Eq.\rref{melnorm}; the sums over $(q, p) \in \D_{h k, m - \ell}$
of the coefficients $P_{m | q p}$ are related as in
Eq.\rref{rel2} to the coefficients $P_{m, h + k, \ell}$
of the polynomial $P_m$, and so
\beq \| \nabla^m \left(\GG(f,\x)\right) \|_{L^2} \leq
\Bigg( \sum_{\stackrel{ h, k \geq 0,~1 \leq h + k \leq m}
{\rule{0mm}{2mm}\scrscr{0 \leq \ell \leq m - h - k}}}
\left( \barray{c} h + k \\ h \farray \right)
P_{m, h + k, \ell}~ \times \feq
$$ \times (1 - {\ell \over m}) U_{m - \ell, h + k, d}
\sigma_{h k \ell}(G, S_{a d} \| f \|_a)
(S_{a d } \| f \|_a)^{h+k-1} \Bigg)~ \| \nabla^m f \|_{L^2} + $$
$$ + \Bigg( \sum_{\stackrel{ h, k \geq 0,~1 \leq h + k \leq m}
{\rule{0mm}{2mm}\scrscr{0 \leq \ell \leq m - h - k}}}
\left( \barray{c} h + k \\ h \farray \right)
P_{m, h + k, \ell} {\ell \over m} U_{m - \ell, h + k, d}
\sigma_{h k \ell}(G, S_{a d} \| f \|_a)
(S_{a d } \| f \|_a)^{h+k-1} \Bigg) \| f \|_{L^2} + $$
$$ + \flat_{m}(G, S_{a d} \| f \|_a )~ \| f \|_{L^2}~. $$
Let us compare this result with the definition \rref{sijlc} of
$\sigma_{j \ell}(G, \rho)$ in terms of the
coefficients $\sigma_{h k \ell}(G, \rho)$,
with $h + k = j$; this gives
\beq \| \nabla^m \left(\GG(f,\x)\right) \|_{L^2} \leq \feq
$$ \leq
\Bigg( \sum_{\stackrel{ 1 \leq j \leq m}
{\rule{0mm}{2mm}\scrscr{0 \leq \ell \leq m - j}}}
P_{m j \ell} (1 - {\ell \over m}) U_{m - \ell, j d}
\sigma_{j \ell}(G, S_{a d} \| f \|_a)
(S_{a d } \| f \|_a)^{j-1} \Bigg)~ \| \nabla^m f \|_{L^2} + $$
$$ + \Bigg( \sum_{\stackrel{ 1 \leq j \leq m}
{\rule{0mm}{2mm}\scrscr{0 \leq \ell \leq m - j}}}
P_{m j \ell} {\ell \over m} U_{m - \ell, j d}
\sigma_{j \ell}(G, S_{a d} \| f \|_a)
(S_{a d } \| f \|_a)^{j-1} \Bigg) \| f \|_{L^2} + $$
$$ + \flat_{m}(G, S_{a d} \| f \|_a )~ \| f \|_{L^2}~. $$
From the expression \rref{form} of the polynomial $P_m$ and
the definitions \rref{betam} \rref{bm} of
$\beta_{m d}$, $b_{m d}$ in terms of $P_m$, we see that
the previous equation means
\beq \| \nabla^m \left(\GG(f,\x)\right) \|_{L^2} \leq~ \mbox{r.h.s. of
Eq.\rref{wget}}~, \feq
as desired.
\vskip 0.1cm \noindent
\textsl{Step 3 (for $n \geq 1$).
It is $\nabla^m \left(G(f,\x)\right) \in L^2(\reali^d, \complessi)$
and \rref{wget} holds for $1 \leq m \leq n$ and $f \in
H^{n}(\reali^d, \complessi) \cap H^a_{r}(\reali^d, \reali)$.}
Let $f \in H^n \cap H^{a}_{r}$, and put
\beq g := G(f, \x) - G(0, \x)~. \label{defg} \feq
This is an $L^2$ function due to Step 1;
the thesis follows if we show that $\nabla^m g$ is $L^2$ for $1 \leq m \leq n$,
with the $L^2$ norm bounded by the r.h.s. of \rref{wget}. As anticipated,
a density argument will be employed to prove this.
\parn
In the sequel, $\ss$ is an index ranging over $\naturali$, and $\vains$
indicates convergence for $\ss \vain +\infty$. By standard
density results, there is a sequence of functions $(f_\ss)$ such that
\beq f_\ss \in H^{n}(\reali^d, \complessi) \cap H^a_{r}(\reali^d, \reali)
\cap C^n(\reali^d, \complessi)~, \qquad \| f_\ss - f \|_{\max(n,a)}
\vains 0~. \feq
Let us put
\beq g_\ss := G(f_\ss, \x) - G(0, \x)~; \feq
then, due to Steps 1 and 2, we have
\beq g_\ss \in L^2(\reali^d, \complessi)~, \qquad
\nabla^m g_\ss \in L^2(\reali^d, \otimes^m \complessi) \quad (1 \leq m \leq n)
\feq
\beq \| g_\ss \|_{L^2} \leq
\flat_{0}(G, S_{a d} \| f_\ss \|_a ) \| f_\ss \|_{L^2} \vains
\flat_{0}(G, S_{a d} \| f \|_a ) \| f \|_{L^2}~,
\label{zkget} \feq
\beq \| \nabla^m g_\ss  \|_{L^2} \leq
\beta_{m d}(G, S_{a d} \| f_\ss \|_a)~ \| \nabla^m f_\ss \|_{L^2}
+ b_{m d}(G, S_{a d} \| f_\ss \|_a)~ \| f_\ss \|_{L^2}
+ \label{wkget} \feq
$$ + \flat_m(G, S_{a d} \| f_\ss \|_a)~ \| f \|_{L^2} \vains
\beta_{m d}(G, S_{a d} \| f \|_a)~ \| \nabla^m f \|_{L^2} + $$
$$ + b_{m d}(G, S_{a d} \| f \|_a)~ \| f \|_{L^2}
+ \flat_m(G, S_{a d} \| f \|_a)~ \| f \|_{L^2} \qquad\qquad (1 \leq
m \leq n)~; $$
the statements on the limits over $\ss$ follow from the continuity of
the mappings
$\flat_0$, $\flat_m$, $b_{m d}$, $\beta_{m d}$, and also imply
\beq \sup_{\ss} \| g_\ss \|_{L^2} < + \infty~, \qquad
\sup_{\ss} \| \nabla^m g_\ss \|_{L^2} < + \infty \quad (1 \leq m \leq n)~.
\label{eqb} \feq
Let us return to the function $g$ of Eq.\rref{defg}, which is also $L^2$
due to Step 1; we claim that
\beq \| g_\ss - g \|_{L^2} \vains 0~. \label{eqc} \feq
In fact $g_\ss - g = G(f_\ss, \x) - G(f, \x)$; to estimate this difference,
we note that, by the Lagrange formula, for any radius $\rho \in [0, r)$ and all
$z, z' \in \Bc(0,\rho)$, $x \in \reali^d$, it is
$$ G(z', x) - G(z, x) = \int_{0}^{1} d t~\Big[~
\partial G((1 - t) z + t z', x)~ (z'- z) +
\overline{\partial} G((1 - t) z + t z', x)~
(\overline{z'} - \overline{z})~ \Big]~, $$
whence
\beq  | G(z', x) - G(z, x) | \leq
\left( \sigma_{1 0 0}(G, \rho) +
\sigma_{0 1 0}(G, \rho) \right) | z' - z |~. \label{app1} \feq
In particular, we choose $\rho$
such that $S_{a d} \| f \|_a < \rho$ and $S_{a d} \| f_\ss \|_a < \rho$
for $\ss$ large; applying Eq.\rref{app1} with $z = f(x)$ and
$z' = f_\ss(x)$, and integrating (the squares of) both sides over $x$,
we get
$$ \| g_\ss - g \|_{L^2} \leq \left( \sigma_{1 0 0}(G, \rho) +
\sigma_{0 1 0}(G, \rho) \right) \| f_\ss - f \|_{L^2} \vains 0~, $$
thus proving \rref{eqc}. The rest of our argument depends on some known facts
on weak convergence ({\footnote{Given a Banach space $X$, a sequence
$(u_\ss)$ with elements in $X$ and another element $u$ of $X$, it is said that
$u_\ss \vains u$ weakly in $X$ if $< \alpha, u_\ss > \vains < \alpha, u >$
for each $\alpha$ in the dual Banach space $X'$. By the uniform boundedness
theorem, this implies $\sup_{\ss} \| u_\ss \|_{X} < +\infty$ and
$\| u \|_{X} \leq \liminf_{\ss} \| u_\ss \|_{X}$, see e.g. \cite{Yos}.
The use of weak convergence in relation to Nemytskij operators is suggested
in \cite{Run}.}}).
Of course, Eq.\rref{eqc} implies
\beq g_\ss \vains g \quad \mbox{weakly in $L^2(\reali^d, \complessi)~.$}
\label{eqx1} \feq
According to a general result on Sobolev spaces (see e.g. \cite{Fri},
Part 1, Sect.6, Lemma 6.2), statements \rref{eqb} \rref{eqx1}
are sufficient to infer that for $1\leq m\leq n$
\beq \nabla^m g \in L^2(\reali^d, \otimes^m \complessi^d), \qquad
\nabla^m g_\ss \vains \nabla^m g \quad
\mbox{weakly in $L^2(\reali^d, \otimes^m \complessi^d)$}~. \feq
So we have (see footnote)
\beq \| \nabla^m g \|_{L^2} \leq \liminf_{\ss} \| \nabla^m g_\ss \|_{L^2} \feq
for $1 \leq m \leq n$; this fact, with Eq.\rref{wkget}, implies
\beq \| \nabla^m g \|_{L^2} \leq~ \mbox{r.h.s. of \rref{wget}} \feq
and the proof is concluded.
\fine
Finally, we are ready to give the
\vskip 0.2cm \noindent
\textbf{Proof of Prop.\ref{nemc}.} From the previous Proposition,
we already know that $G(f, \x)$ is $H^n$; let us evaluate
\beq \| G(f, \x) - G(0, \x) \|_{n} =
\sqrt{\sum_{m=0}^n \left( \barray{c} n \\ m \farray \right)
\| \nabla^m (G(f, \x)) - \nabla^m (G(0, \x)) \|^{2}_{L^2}}~. \label{nonagf} \feq
Prop.\ref{dmg} tells us that
\beq \| \nabla^m \left(G(f,\x)\right) - \nabla^m \left(G(0,\x)\right)
\|_{L^2} \leq X_m + Y_m
\qquad (m=0,...,n) \feq
$$ X_0 := \flat_{0}(G, S_{a d} \| f \|_a )~\| f \|_{L^2}, \quad X_m :=
\beta_{m d}(G, S_{a d} \| f \|_a)~\| \nabla^m f \|_{L^2}
\quad (1 \leq m \leq n)~, $$
$$ Y_0 := 0,
\quad Y_m := \left(b_{m d}(G, S_{a d} \| f \|_a)
+ \flat_m(G, S_{a d} \| f \|_a) \right)~ \| f \|_{L^2}~~(1 \leq m \leq n). $$
Inserting this into Eq.\rref{nonagf} we get
\beq \| G(f, \x) - G(0, \x) \|_{n} \leq
\sqrt{\sum_{m=0}^n \left( \barray{c} n \\ m \farray \right)
(X_m + Y_m)^2 } \leq  \label{gg1} \feq
$$ \leq
\sqrt{\sum_{m=0}^n \left( \barray{c} n \\ m \farray \right) X_m^2 } +
\sqrt{\sum_{m=0}^n \left( \barray{c} n \\ m \farray \right) Y_m^2 }~.
$$
On the other hand, the definitions of $X_m$, $Y_m$ imply
\beq\sqrt{\sum_{m=0}^n \left( \barray{c} n \\ m \farray \right) X_m^2 }  \leq
\max \left\{ \flat_{0}(G, S_{a d} \| f \|_a),~\beta_{m d}(G, S_{a d} \| f \|_a)~
(m=1,...,n) \right\}~\| f \|_n~; \label{xn1} \feq
\beq
\sqrt{\sum_{m=0}^n \left( \barray{c} n \\ m \farray \right) Y_m^2 } =
\sqrt{\sum_{m=1}^n \left( \barray{c} n \\ m \farray \right)
\left( b_{m d}(G, S_{a d} \| f \|_a) +
\flat_m(G, S_{a d} \| f \|_a) \right)^2}~\| f \|_{L^2}~;
\label{yn1} \feq
inserting \rref{xn1} \rref{yn1} into \rref{gg1}, we obtain the bound
(\ref{disug}-\ref{cnd}) on $\| G(f, \x) - G(0, \x) \|_n$;
the weaker bound \rref{weak} follows
trivially, since $\| f \|_{L^2} \leq \| f \|_n$. \fine
\vskip 0.2cm\noindent
\textbf{Comments and acknowledgments.} The
quantitative estimates in the present work
are extended with little effort to
Nemytskij operators on spaces $H^n(\reali^d, \reali^\delta$  ~$\mbox{or}~
\complessi^{\delta})$, where $\delta$ is an arbitrary integer; here,
this has not been done to keep notations as simple as possible.
We again acknowledge
W. Beckner, L. Colzani, G. Meloni and S. Paveri Fontana
for useful indications. \parn
This work was partly supported by the GNFM
of Istituto Nazionale di Alta Matematica and by MURST.
\vskip 0.2cm \noindent
\appendix
\section{Appendix. Some facts about tensors.}
\label{appe}
\textbf{Identities for permutation and symmetrization operators.}
For each $\ell$, we have defined the permutation
operators $\PP_{\sigma}$ and the symmetrization operator
$\SS$, sending $\otimes^\ell \complessi^d$ into itself
(Eq.s \rref{perm} \rref{simm}). For all $\ell$
and all permutations $\sigma, \sigma'$
$\in \boma{\ell!}$, it is readily checked that
\beq \PP_{\sigma} \PP_{\sigma'} = \PP_{\sigma \circ \sigma'}~,
\feq
\beq \SS \PP_{\sigma} = \PP_{\sigma} \SS = \SS~.
\label{use2} \feq
From here, with elementary manipulations one infers the relations
\beq \SS((\SS T) \otimes U) = \SS (T \otimes (\SS U)) = \SS(T \otimes U)
\label{man} \feq
for all tensors $T \in \otimes^\ell
\complessi^d$, $U \in \otimes^m \complessi^d$
(see e.g. \cite{Kos}, Chapter 4, in the proof of Prop. 5.7).
\vskip 0.1cm \noindent
\textbf{Proofs of the commutative
and associative properties \rref{assoc} of the symmetrized
tensor product.} The arguments are standard \cite{Kos},
and reported only for completeness. Let
$T \in \otimes^\ell \complessi^d$, $U \in \otimes^m \complessi^d$.
It is evident that $U \otimes T = \PP_{\zeta}(T \otimes U)$,
where $\zeta \in \boma{(\ell + m)!}$ is the permutation
such that $\zeta(1) := m+1$, ..., $\zeta(\ell) := m + \ell$,
$\zeta(\ell + 1) := 1$, ..., $\zeta(\ell + m) := m$. Therefore
$U \vee T = \SS \PP_{\zeta}(T \otimes U)$ $= \SS(T \otimes U)$
$= T \vee U$ ($\PP_{\zeta}$ disappears due to Eq. \rref{use2}).
\parn
Now, let us consider a third tensor
$V \in \otimes^p \complessi^d$;
then $(T \vee U) \vee V$ $ = \SS(\SS(T \otimes U) \otimes V)$
$ = \SS((T \otimes U) \otimes V))$ (the first passage holds
by the definition of $\vee$, the second follows from
\rref{man}). In a similar way one finds $T \vee (U \vee V) $
$= \SS(T \otimes (U \otimes V))$; now, the associativity
of $\otimes$ yields the equality $(T \vee U) \vee V =
T \vee (U \vee V)$. \parn
\vskip 0.1 cm \noindent
\textbf{Proof of Eq.\rref{notv} for the tensor norms.}
Consider tensors $T \in \otimes^\ell \complessi^d$,
$U \in \otimes ^m \complessi^d$. The equality
$| T \otimes U |^2$  $= | T |^2 | U |^2$ is checked
in an elementary way, explicitating the definitions of
$| ~|$ and $\otimes$. Furthermore,
$$ | T \vee U | = | \SS(T \otimes U) |
\leq {1 \over (\ell + m)!} \sum_{\sigma \in \boma{(\ell+ m) !}}
| \PP_{\sigma} (T \otimes U) |~, $$
the last passage following
from the definition of $\SS$ and the triangular inequality for
$|~|$. On the other hand, it is clear from
the very definition that each operator $\PP_{\sigma}$
is isometric w.r.t. the norm $|~|$; so
$| \PP_{\sigma} (T \otimes U) |$  $= | T \otimes U |$  $= | T | | U |$
for all $\sigma \in \boma{(\ell + m)!}$, and
this implies $| T \vee U | \leq | T | | U |$, as desired.
\vskip 0.1cm \noindent
\textbf{Identities for the derivative $\nabla$
and the operators $\boma{\PP_{\sigma}}$, $\boma{\SS}$.}
We regard all these operators as acting on the
$C^1$ tensor fields of some
definite order. For each order $\ell $ and each
permutation $\sigma \in \boma{\ell !}$, it is
\beq \nabla \PP_{\sigma} = \PP_{\widetilde{\sigma}} \nabla
\label{use1} \feq
where $\widetilde{\sigma} \in \boma{(\ell + 1)!}$ is
defined by $\widetilde{\sigma}(i) := \sigma(i)$ for
$1 \leq i \leq \ell$,~$\widetilde{\sigma}(\ell + 1) :=
\ell+1$; this fact is straightforwardly checked using the
definition \rref{nab} of $\nabla$. \parn
As a second fact, on tensor fields of any order $\ell$ we
have the identity
\beq \SS \nabla \SS = \SS \nabla~; \label{ident} \feq
this follows writing $\SS \nabla \SS$
$= (1/\ell!) \sum_{\sigma \in \boma{\ell!}}$
${\SS \nabla \PP_{\sigma}}$, and using Eq.s \rref{use1}
\rref{use2}.
\vskip 0.1cm \noindent
\textbf{Leibnitz rule for $\boma{\nabla}$ and $\boma{\otimes}$.}
For all $C^1$ tensor fields $T : \reali^d
\vain \otimes^\ell \complessi^d$ and
$U : \reali^d \vain \otimes^m \complessi^d$, it is
\beq \nabla(T \otimes U) = \PP_{\eta} (\nabla T \otimes U) +
T \otimes \nabla U~, \label{compa} \feq
where $\eta \in \boma{(\ell+m+1)!}$ is defined by
$\eta(1) := 1$,...,$\eta(\ell) := \ell$, $\eta(\ell+1) := \ell+m+1$,
..., $\eta(\ell+2) := \ell+1$, $\eta(\ell+3) := \ell+2$, ...,
$\eta(\ell+m+1) := \ell + m$.
To prove this we note that, in terms of components, we have
$$ \left(\nabla(T \otimes U)\right)_{\lambda_1 ... \lambda_{\ell+m+1}} =
\partial_{\lambda_{\ell+m+1}} \left( T_{\lambda_1 ... \lambda_\ell}
U_{\lambda_{\ell+1} ... \lambda_{\ell+m}} \right) = $$
$$ = \left( \partial_{\lambda_{\ell+m+1}} T_{\lambda_1 ... \lambda_\ell}
\right) U_{\lambda_{\ell+1}...\lambda_{\ell+m}} +
T_{\lambda_1...\lambda_\ell}~
\left(\partial_{\ell+m+1} U_{\lambda_{\ell+1}...\lambda_{\ell+m}}
\right)~; $$
Eq.\rref{compa} is just a compact formulation of this
result.
\vskip 0.1cm \noindent
\textbf{"Symmetrized" Leibnitz rule for $\boma{\nabla_S}$ and $\boma{\vee}$.}
This is given by Eq. \rref{lei2}, here reported:
$$ \nabla_{\SS}(T \vee U) =
(\nabla_{\SS} T) \vee U + T \vee (\nabla_{\SS} U)~$$
for all $C^1$ tensor fields $T, U$ of arbitrary orders
$\ell, m$. Eq. \rref{lei2} is proved by the chain of relations
$$ \nabla_{\SS}(T \vee U) \stackrel{\mbox{\footnotesize{def.s of}
$\nabla_{\SS}$, $\vee$}}{=} \SS \nabla \SS(T \otimes U)
\stackrel{\mbox{\footnotesize{\rref{ident}}} }{=} \SS \nabla(T \otimes U) = $$
$$ \stackrel{\mbox{\footnotesize{\rref{compa}}}}{=} \SS \PP_{\eta}(\nabla T
\otimes U) + \SS(T \otimes \nabla U)
\stackrel{\mbox{\footnotesize{\rref{use2}}}}{=}
\SS(\nabla T
\otimes U) + \SS(T \otimes \nabla U) = $$
$$ \stackrel{\mbox{\footnotesize{\rref{man}}}}{=}
\SS(\SS \nabla T
\otimes U) + \SS(T \otimes \SS \nabla U)
\stackrel{\mbox{\footnotesize{def.s of}
$\nabla_{\SS}$, $\vee$}}{=} \nabla_{\SS} T \vee U + T \vee \nabla_{\SS} U~. $$


\vskip 0.2cm \noindent
\begin{thebibliography}{99}
\bibitem{Ada2} D.R. Adams, \textsl{On the existence of capacitary strong type
estimates in $\reali^n$}, Ark. Math. \textbf{14}, 125-140 (1976);
D.R. Adams, M. Frazier: \textsl{BMO and smooth truncation in Sobolev spaces},
Studia Math. \textbf{39}, 241-260 (1988),
\textsl{Composition operators on potential spaces},
Proc. AMS \textbf{114}, 155-165 (1992)
\bibitem{Smi} N. Aronszajn, K.T. Smith, \textsl{Theory of Bessel
potentials. I}, Ann. Inst. Fourier \textbf{11}, 385-475 (1961)
\bibitem{Bec} W. Beckner, \textsl{Inequalities in Fourier
Analysis}, Ann. of Math. \textbf{102}, 159-182 (1975)
\bibitem{Bou} G. Bourdeau: \textsl{Le calcul fonctionnel dans les espaces
de Sobolev}, Invent. math. \textbf{104}, 435-446 (1991),
\textsl{The functional calculus in Sobolev spaces}, Proc. Conf.: Function
spaces, differential operators and nonlinear analysis, Teubner-Text Math.
\textbf{133}, 127-142, Teubner, Leipzig (1993)
\bibitem{Com} L. Comtet, \textsl{Advanced Combinatorics}, Reidel, Dordrecht
(1974)
\bibitem{Fri} A. Friedman, \textsl{Partial differential equations},
Holt, Rinehart and Winston, New York (1969)
\bibitem{Gar} J. Garcia-Cuerva, J.L. Rubio de Francia, \textsl{Weighted
norm inequalities and related topics}, Mathematics Studies
\textbf{116}, North-Holland, Amsterdam (1985)
\bibitem{Gzy} H. Gzyl, \textsl{Multidimensional extension of Fa\`a di Bruno's
formula}, J. Math. Anal. Appl. \textbf{116}, 450-455 (1986)
\bibitem{Ham} R.S. Hamilton, \textsl{The inverse function theorem of Nash
and Moser}, Bull. Amer. Math. Soc. (N.S.) \textbf{7}, 65-222 (1982)
\bibitem{Kac} V. Kac, \textsl{Infinite dimensional Lie algebras}, Cambridge
Univ. Press, Cambridge (1990)
\bibitem{Kos} A. I. Kostrikin, Y.I. Manin, \textsl{Linear algebra and
geometry}, Gordon and Breach, New York (1989)
\bibitem{Lie2} E. H. Lieb, M. Loss, \textsl{Analysis}, Graduate
Studies in Mathematics \textbf{14}, American Mathematical Society
(1997)
\bibitem{Maz} V.G. Mazjia, \textsl{Sobolev spaces}, Springer, Berlin (1985)
\bibitem{Mel} G. Meloni, private communication
\bibitem{Mis} R. L. Mishkov, \textsl{Generalization of the formula of Fa\`a
di Bruno for a composite function with a vector argument}, Int. J. Math.
Math. Sci. \textbf{24}, 481-491 (2000)
\bibitem{quad} C. Morosi, L. Pizzocchero, \textsl{Semilinear evolution
equations in Fr\'echet spaces. Abstract theory}, Quaderno 23/1999,
Dipartimento di Matematica, Universit\`a di Milano
\bibitem{mp1} C. Morosi, L. Pizzocchero, \textsl{On the constants for
some Sobolev imbeddings}, J. Inequal. Appl. \textbf{6}, 665-679
(2001)
\bibitem{mp2} C. Morosi, L. Pizzocchero, \textsl{On the constants in
some inequalities for the Sobolev norms and pointwise product}, J.
Inequal. Appl. \textbf{7}, 421-452 (2002)
\bibitem{Run} T. Runst, W. Sickel, \textsl{Sobolev spaces of
fractional order, Nemytskij operators and nonlinear partial
differential equations}, de Gruyter, Berlin (1996)
\bibitem{Ster} S. Sternberg, \textsl{Lectures on differential geometry},
Chelsea, New York (1983)
\bibitem{Yos} K. Yosida, \textsl{Functional Analysis}, Springer, Berlin (1965)
\end{thebibliography}
\end{document}